\documentclass[12pt]{article}
\pdfoutput=1

\usepackage{draft}

\newtheorem{thm}{Theorem}[section]
\newtheorem{defi}{Definition}[section]

\begin{document}

\begin{titlepage}

\preprint{CALT-TH-2020-018}

\begin{center}

\hfill \\
\hfill \\
\vskip 1cm

\title{
The $F$-Symbols for Transparent Haagerup-Izumi \\ Fusion Categories with $G = \mathbb{Z}_{2n+1}$
}

\author{Tzu-Chen Huang, Ying-Hsuan Lin}

\address{
Walter Burke Institute for Theoretical Physics,
\\ 
California Institute of Technology, Pasadena, CA 91125, USA
}

\email{jimmy@caltech.edu, yinhslin@gmail.com}

\end{center}

\vfill

\begin{abstract}
A fusion category is called transparent if the associator involving any invertible object is the identity map.  For the Haagerup-Izumi fusion rings with $G = \mathbb{Z}_{2n+1}$ (the $\mathbb{Z}_3$ case is the Haagerup fusion ring with six simple objects), the transparent ansatz reduces the number of independent $F$-symbols from order $\mathcal{O}(n^6)$ to $\mathcal{O}(n^2)$, rendering the pentagon identity practically solvable.  Transparent Haagerup-Izumi fusion categories are thereby constructively classified up to $G = \mathbb{Z}_9$, recovering all known Haagerup-Izumi fusion categories to this order, and producing new ones.  Transparent Haagerup-Izumi fusion categories additionally satisfying $S_4$ tetrahedral invariance are further classified up to $G = \mathbb{Z}_{15}$, and the explicit $F$-symbols for the unitary ones, including the Haagerup $\mathcal{H}_3$ fusion category, are compactly presented.  The $F$-symbols for the Haagerup $\mathcal{H}_2$ fusion category are also presented. Going beyond, the transparent ansatz offers a viable course towards constructing novel fusion categories for new fusion rings.
\end{abstract}

\vfill

\end{titlepage}

\tableofcontents

\section{Introduction}

Subfactors \cite{Jones:1983kv,jones1997introduction} and fusion categories \cite{Etingof:aa,etingof2016tensor} provide the mathematical framework underlying various physical objects in quantum field theory, including anyons in (2+1)$d$ Chern-Simons theory \cite{Witten:1988hf,jones1990neumann} and topological defect lines in (1+1)$d$ quantum field theory \cite{Bhardwaj:2017xup,Tachikawa:2017gyf,Chang:2018iay}.  Fusion categories with invertible objects encapsulate the notion of symmetries and 't Hooft anomalies in quantum field theory, and those with non-invertible objects generalize such a notion\cite{Chang:2018iay,tHooft:1980xss,Kapustin:2014zva,Thorngren:2019iar}.  Due to Ocneanu rigidity \cite{wassermann2000quantum,Etingof:aa}, a fusion category is an invariant under renormalization group flows connecting short and long distance physics.  This generalization of the 't Hooft anomaly matching condition has shed new light on the phases of quantum field theory.

Subfactor theory has an inherent categorical structure \cite{muger2003subfactors}, and has been a productive factory of fusion categories.  Subfactors with Jones indices less than 4 have been classified by Ocneanu \cite{ocneanu1988quantized} and extended to 4 by Popa \cite{popa1994classification}.  Haagerup \cite{haagerup1994principal} searched for subfactors with Jones indices {\it a little bit beyond 4}, and together with Asaeda \cite{asaeda1999exotic} constructed one with Jones index ${5+\sqrt{13}\over2}$, the smallest above 4.  In \cite{izumi2001structure}, Izumi generalized the Haagerup fusion ring to a family of fusion rings labeled by a finite abelian group $G$, and explicitly constructed the subfactors for $G = \bZ_3, \, \bZ_5$.  The constructive classification of subfactors for $|G|$ odd was achieved up to $|G| = 19$ by Evans and Gannon \cite{Evans:2010yr} (up to $|G| = 9$ with exact expressions and the rest with numerical estimates), and that of subfactors with $G = \bZ_4, ~ \bZ_2 \times \bZ_2, ~ \bZ_4 \times \bZ_2, ~ \bZ_6, ~ \bZ_8, ~ \bZ_{10}$ by Grossman, Izumi, and Snyder \cite{grossman2015drinfeld,grossman2018asaeda,Izumi_2018,grossman2019infinite}.

The fundamental data underlying a fusion category are the $F$-symbols, which are solutions to the pentagon identity.  Some (almost) equivalent notions exist:  associators, quantum $6j$-symbols and crossing kernels.  They underlie the Turaev-Viro theory \cite{Turaev:1992hq,Turaev:1994xb}, the Levin-Wen string-net models \cite{Levin:2004mi}, and large classes of statistical models (see \cite{Aasen:2020jwb} and references within) as well as the associated anyon chains \cite{Feiguin:2006ydp}.  In \cite{Chang:2018iay}, one of the present authors showed how the $F$-symbols strongly constrain (1+1)$d$ (fully extended) topological quantum field theories \cite{Moore:2006dw,Davydov:2011kb}, which are endpoints of symmetry-preserving renormalization group flows; in many cases, given the $F$-symbols, the full field theory data could be completely determined by bootstrap.

In this paper, unitary and non-unitary Haagerup-Izumi fusion categories with $G = \bZ_{2n+1}$ are constructed up to $G = \bZ_{15}$ by computing Groebner bases for the pentagon identity.  The notion of a {\it transparent} fusion category is introduced in Definition~\ref{Transparent}, from which various consequential graph equivalences and $F$-symbol relations are derived to reduce the number of independent $F$-symbols from $\mathcal{O}(n^6)$ to $\mathcal{O}(n^2)$ and render the pentagon identity practically solvable.  These relations are summarized into a system of constraints in Definition~\ref{T}, and the solutions to the pentagon identity under said constraints provide a classification of $F$-symbols for transparent Haagerup-Izumi fusion categories.  The results of this classification are stated in Theorems~\ref{Main} and~\ref{Main2}.

Some remarks on the comparison of the present results with the existing literature are in order.  As mentioned above, the datum equivalent to the $F$-symbols for several unitary Haagerup-Izumi fusion categories were obtained by Izumi \cite{izumi2001structure}, Evans and Gannon \cite{Evans:2010yr}, and Grossman and Snyder \cite{Grossman_2012} using Cuntz algebra techniques; such constructions were further generalized by Evans and Gannon \cite{Evans:2015zga} to fusion categories that need not be unitary.  More recently, the $F$-symbols for all fusion categories realizing the Haagerup fusion ring ($G = \bZ_3$) with six simple objects have been computed using the pentagon approach by Titsworth \cite{Titsworth}, and for the special case of the Haagerup $\cH_3$ fusion category (in the nomenclature of Grossman and Snyder \cite{Grossman_2012}) independently by Osborne, Stiegemann and Wolf \cite{osborne2019fsymbols}.

The novelty of this paper is twofold.  First, it offers the direct pentagon construction for Haagerup-Izumi fusion categories beyond the Haagerup case ($G = \bZ_3$); in particular, the Haagerup-Izumi fusion categories classified in Theorem~\ref{Main2} have not appeared in the literature beyond $G = \bZ_5$.  Second, the special transparent gauge adopted in this paper---in which all $F$-symbols involving at least one external invertible object take value one---not only makes the independent $F$-symbols directly comparable to the Cuntz algebra datum of Izumi \cite{izumi2001structure}, Evans and Gannon \cite{Evans:2010yr,Evans:2015zga}, and Grossman and Snyder \cite{Grossman_2012}, but also makes the $F$-symbols automatically tetrahedral-symmetric ($A_4$ or $S_4$ tetrahedral-invariant in the language of this paper), and unitary for pseudo-unitary fusion categories.\footnote{In \cite{Titsworth,osborne2019fsymbols}, the $F$-symbols for the Haagerup fusion categories with six simple objects were presented in non-transparent gauges that do not enjoy tetrahedral symmetry.  The present authors used the Mathematica package provided by Titsworth \cite{Titsworth} to check that the $F$-symbols in the present paper are indeed gauge-equivalent to his.  The authors also thank Yuji Tachikawa for explicitly checking that the four sets of $F$-symbols in \cite{osborne2019fsymbols} are all gauge-equivalent, and also gauge-equivalent to those presented in this paper.
}
In physical applications, such a gauge satisfies the assumptions of various theoretical constructions---the Levin-Wen string-net models \cite{Levin:2004mi}, large classes of statistical models (see \cite{Aasen:2020jwb} and references within) and the associated anyon chains \cite{Feiguin:2006ydp}---and allows the more effective exploitation of the $G = \bZ_{2n+1}$ symmetry.  Of course, for a given fusion ring, there may exist non-transparent fusion categories that elude the present approach.  However, none of the Haagerup-Izumi fusion categories up to $G = \bZ_9$ known in the literature \cite{asaeda1999exotic,izumi2001structure,Evans:2010yr,grossman2015drinfeld,grossman2018asaeda,Izumi_2018,grossman2019infinite} was found to be non-transparent!

The outline of this paper is as follows.  Section~\ref{Sec:Prelim} reviews the string diagram calculus, the $F$-symbols, and their relation to the tetrahedra.  Section~\ref{Sec:Transparent} defines the notion of a transparent fusion category, and derives various consequences including invariance relations for the $F$-symbols.  Section~\ref{Sec:HI} introduces the Haagerup-Izumi fusion rings, and formulates a set of constraints on $F$-symbols that must be satisfied for transparent Haagerup-Izumi fusion categories.  Section~\ref{Sec:F} states the classification of solutions to the pentagon identity under the said constraints, and presents the explicit $F$-symbols for unitary Haagerup-Izumi fusion categories with $S_4$ tetrahedral invariance, as well as for the Haagerup $\cH_2$ fusion category.  Section~\ref{Sec:Conclusions} ends with some concluding remarks.

{\it Note:} The authors first obtained the $F$-symbols for the Haagerup fusion categories with six simple objects from Titsworth \cite{Titsworth}.  By performing gauge transformations on his solution, a gauge manifesting the transparent property was found.  This observation led the present authors to postulate that transparent fusion categories also exist for the subsequent Haagerup-Izumi fusion rings with $G = \bZ_{2n+1}$.

\section{Preliminaries}
\label{Sec:Prelim}

A classic introduction to fusion categories can be found in \cite{Etingof:aa,etingof2016tensor}.  The type of fusion categories considered in this paper are pivotal fusion categories over ground field $k = \bC$.\footnote{For such categories, a physical formulation in the context of topological defect lines in (1+1)$d$ quantum field theory can be found in \cite{Chang:2018iay} (see also \cite{Bhardwaj:2017xup,Tachikawa:2017gyf}).  
}
The notation for string diagrams is as follows.  Each object $\cL$ is represented by an oriented string that is equivalent to its dual $\ocL$ with the opposite orientation,
\[
\begin{gathered}
\begin{tikzpicture}[scale=1]
\draw [line,->-=.6] (0,0) node [below] {} -- (0,1) node [above] {$\cL$};
\end{tikzpicture}
\end{gathered}
\quad = \quad
\begin{gathered}
\begin{tikzpicture}[scale=1]
\draw [line,->-=.6] (0,1) node [above] {$\ocL$} -- (0,0) node [below=1mm] {};
\end{tikzpicture}
\end{gathered} \, .
\]
The basic building block for string diagrams is a trivalent vertex with three open edges
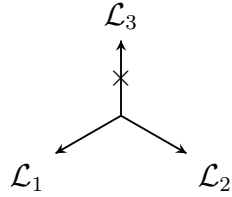
\[
\begin{gathered}
\begin{tikzpicture}[scale=1]
\draw [line,->-=1] (0,0) -- (0,.5) node {$\times$} -- (0,1) node [above] {$\cL_3$};
\draw [line,->-=1] (0,0) -- (-.87,-.5) node [below left] {$\cL_1$};
\draw [line,->-=1] (0,0) -- (.87,-.5) node [below right] {$\cL_2$};
\end{tikzpicture}
\end{gathered}
\]
with $\times$ specifying the ordering of edges.  It represents the vector space of morphisms
\[
V_{\cL_1, \cL_2, \cL_3} \equiv \hom(\ocL_2 \otimes \ocL_1, \cL_3) \in \bC^{N_{\ocL_2, \ocL_1}^{\cL_3}} \, ,
\]
where $N_{\ocL_2, \ocL_1}^{\cL_3}$ is the fusion coefficient, the multiplicity of $\cL_3$ in $\ocL_2 \otimes \ocL_1$.  A change of basis at this vertex is a {\it gauge} transformation $g_{\cL_1, \cL_2, \cL_3} \in GL(N_{\cL_1, \cL_2}^{\cL_3}, \bC)$.  To simplify the discussion, it is assumed in the following that the fusion algebra is multiplicity-free, {\it i.e.} all nonzero fusion coefficients are one, hence every nontrivial gauge factor $g_{\cL_1, \cL_2, \cL_3}$ is a complex scalar.

For a trivalent vertex involving at least one unit object, the ordering of edges is irrelevant, and the marking $\times$ can be dropped.  Furthermore, by choosing the unitors and counitors to be identity morphisms, the unit object $\cI$ can be removed or added at will,
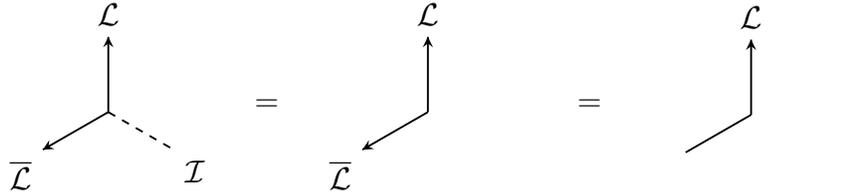
\[
\begin{gathered}
\begin{tikzpicture}[scale=1]
\draw [line,->-=1] (0,0) -- (0,1) node [above] {$\cL$};
\draw [line,->-=1] (0,0) -- (-.87,-.5) node [below left] {$\ocL$};
\draw [line,dashed] (0,0) -- (.87,-.5) node [below right] {$\cI$};
\end{tikzpicture}
\end{gathered}
\quad = \quad
\begin{gathered}
\begin{tikzpicture}[scale=1]
\draw [line,->-=1] (0,0) -- (0,1) node [above] {$\cL$};
\draw [line,->-=1] (0,0) -- (-.87,-.5) node [below left] {$\ocL$};
\draw [line,draw=none] (0,0) -- (.87,-.5) node [below right=3mm] {};
\end{tikzpicture}
\end{gathered}
\quad = \quad
\begin{gathered}
\begin{tikzpicture}[scale=1]
\draw [line,->-=1] (0,0) -- (0,1) node [above] {$\cL$};
\draw [line] (0,0) -- (-.87,-.5) node [below left=3mm] {~};
\draw [line,draw=none] (0,0) -- (.87,-.5) node [below right=3mm] {};
\end{tikzpicture}
\end{gathered} \, .
\]

For a string diagram composed of two trivalent vertices
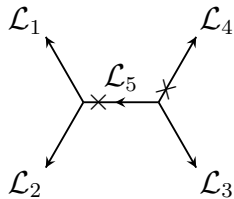
\[
\begin{gathered}
\begin{tikzpicture}[scale=1]
\draw [line,-<-=.56] (-1,0) -- (-.8,0) node {$\times$} -- (-.5,0) node [above] {$\cL_5$} -- (0,0);
\draw [line,->-=1] (-1,0) -- (-1.5,.87) node [above left=-3pt] {$\cL_1$};
\draw [line,->-=1] (-1,0) -- (-1.5,-.87) node [below left=-3pt] {$\cL_2$};
\draw [line,->-=1] (0,0) -- (.5,-.87) node [below right=-3pt] {$\cL_3$};
\draw [line,->-=1] (0,0) -- (.1,.175) node {\rotatebox[origin=c]{60}{$\times$}} -- (.5,.87) node [above right=-3pt] {$\cL_4$};
\end{tikzpicture}
\end{gathered}
\]
the gauge freedom is $g_{\cL_1, \cL_2, \ocL_5} \, g_{\cL_5, \cL_3, \cL_4}$.  It is related by an {\it $F$-move} to a sum of string diagrams in a different configuration,
\ie
\label{F}
\begin{gathered}
\begin{tikzpicture}[scale=1]
\draw [line,-<-=.56] (-1,0) -- (-.8,0) node {$\times$} -- (-.5,0) node [above] {$\cL_5$} -- (0,0);
\draw [line,->-=1] (-1,0) -- (-1.5,.87) node [above left=-3pt] {$\cL_1$};
\draw [line,->-=1] (-1,0) -- (-1.5,-.87) node [below left=-3pt] {$\cL_2$};
\draw [line,->-=1] (0,0) -- (.5,-.87) node [below right=-3pt] {$\cL_3$};
\draw [line,->-=1] (0,0) -- (.1,.175) node {\rotatebox[origin=c]{60}{$\times$}} -- (.5,.87) node [above right=-3pt] {$\cL_4$};
\end{tikzpicture}
\end{gathered}
\quad = \quad
\sum_{\cL_6} (F^{\cL_1, \cL_2, \cL_3}_{\ocL_4})_{\cL_5, \cL_6}
\begin{gathered}
\begin{tikzpicture}[scale=1]
\draw [line,-<-=.56] (0,-1) -- (0,-.8) node {$\times$} -- (0,-.5) node [right] {$\cL_6$} -- (0,0);
\draw [line,->-=1] (0,-1) -- (.87,-1.5) node [below right=-3pt] {$\cL_3$};
\draw [line,->-=1] (0,-1) -- (-.87,-1.5) node [below left=-3pt] {$\cL_2$};
\draw [line,->-=1] (0,0) -- (-.87,.5) node [above left=-3pt] {$\cL_1$};
\draw [line,->-=1] (0,0) -- (.175,.1) node {\rotatebox[origin=c]{30}{$\times$}} -- (.87,.5) node [above right=-3pt] {$\cL_4$};
\end{tikzpicture}
\end{gathered} \, ,
\fe
where $(F^{\cL_1, \cL_2, \cL_3}_{\ocL_4})_{\cL_5, \cL_6}$ are the {\it $F$-symbols}.  The gauge factor for an $F$-symbol is
\[
{g_{\cL_1, \cL_2, \ocL_5} \, g_{\cL_5, \cL_3, \cL_4} \over g_{\cL_2, \cL_3, \ocL_6} \, g_{\cL_1, \cL_6, \cL_4}} \, .
\]

The $F$-symbols must satisfy a consistency condition that is the equivalence of the two different sequences of $F$-moves
\ien
&
\begin{gathered}
\begin{tikzpicture}[scale=.5]
\draw [line,->-=1] (0,0) -- (0,-.282) node {$\times$} -- (0,-1.41);
\draw [line,->-=.56] (0,0) -- (.8,.8) node {$+$} -- (1,1);
\draw [line,->-=.56] (1,1) -- (1.8,1.8) node {$+$} -- (2,2);
\draw [line,->-=1] (2,2) -- (3,3);
\draw [line,->-=1] (0,0) -- (-3,3);
\draw [line,->-=1] (1,1) -- (-1,3);
\draw [line,->-=1] (2,2) -- (1,3);
\end{tikzpicture}
\end{gathered}
\quad \longrightarrow \quad
\begin{gathered}
\begin{tikzpicture}[scale=.5]
\draw [line,->-=1] (0,0) -- (0,-.282) node {$\times$} -- (0,-1.41);
\draw [line,->-=.56] (0,0) -- (1.8,1.8) node {$+$} -- (2,2);
\draw [line,->-=1] (2,2) -- (3,3);
\draw [line,->-=.56] (0,0) -- (-1,1) node [below left=-3pt] {$\cL_8$} -- (-1.8,1.8) node {$+$} -- (-2,2);
\draw [line,->-=1] (-2,2) -- (-3,3);
\draw [line,->-=1] (-2,2) -- (-1,3);
\draw [line,->-=1] (2,2) -- (1,3);
\end{tikzpicture}
\end{gathered}
\quad \longrightarrow \quad
\begin{gathered}
\begin{tikzpicture}[xscale=-.5, yscale=.5]
\draw [line,->-=1] (0,0) -- (0,-.282) node {$\times$} -- (0,-1.41);
\draw [line,->-=.56] (0,0) -- (.8,.8) node {$+$} -- (1,1);
\draw [line,->-=.56] (1,1) -- (1.8,1.8) node {$+$} -- (2,2);
\draw [line,->-=1] (2,2) -- (3,3);
\draw [line,->-=1] (0,0) -- (-3,3);
\draw [line,->-=1] (1,1) -- (-1,3);
\draw [line,->-=1] (2,2) -- (1,3);
\end{tikzpicture}
\end{gathered}
\\
\\
& \hspace{.75in} \searrow \hspace{3in} \nearrow
\\
\\
& \hspace{1in}
\begin{gathered}
\begin{tikzpicture}[scale=.5]
\draw [line,->-=1] (0,0) -- (0,-.282) node {$\times$} -- (0,-1.41);
\draw [line,->-=.56] (0,0) -- (.8,.8) node {$+$} -- (1,1);
\draw [line,->-=1] (1,1) -- (3,3);
\draw [line,->-=1] (0,0) -- (-3,3);
\draw [line,->-=.56] (1,1) -- (.5,1.5) node [above right=-3pt] {$\cL$} -- (.2,1.8) node {$+$} -- (0,2);
\draw [line,->-=1] (0,2) -- (-1,3);
\draw [line,->-=1] (0,2) -- (1,3);
\end{tikzpicture}
\end{gathered}
\quad \longrightarrow \quad
\begin{gathered}
\begin{tikzpicture}[scale=.5]
\draw [line,->-=1] (0,0) -- (0,-.282) node {$\times$} -- (0,-1.41);
\draw [line,->-=1] (0,0) -- (3,3);
\draw [line,->-=.56] (0,0) -- (-.5,.5) node [below left=-3pt] {$\cL_9$} -- (-.8,.8) node {$+$} -- (-1,1);
\draw [line,->-=1] (-1,1) -- (-3,3);
\draw [line,->-=.56] (-1,1) -- (-.5,1.5) node [above left=-3pt] {$\cL$} -- (-.2,1.8) node {$+$} -- (0,2);
\draw [line,->-=1] (0,2) -- (-1,3);
\draw [line,->-=1] (0,2) -- (1,3);
\end{tikzpicture}
\end{gathered}
\fen
both resulting in
\[
\begin{gathered}
\begin{tikzpicture}[scale=.5]
\draw [line,->-=1] (0,0) -- (0,-.282) node {$\times$} -- (0,-1.41) node [below] {$\cL_5$};
\draw [line,->-=.56] (0,0) -- (.5,.5) node [below right=-3pt] {$\cL_7$} -- (.8,.8) node {$+$} -- (1,1);
\draw [line,->-=.56] (1,1) -- (1.5,1.5) node [below right=-3pt] {$\cL_6$} -- (1.8,1.8) node {$+$} -- (2,2);
\draw [line,->-=1] (2,2) -- (3,3) node [above right=-3pt] {$\cL_1$};
\draw [line,->-=1] (0,0) -- (-3,3) node [above left=-3pt] {$\cL_4$};
\draw [line,->-=1] (1,1) -- (-1,3) node [above left=-3pt] {$\cL_3$};
\draw [line,->-=1] (2,2) -- (1,3) node [above left=-3pt] {$\cL_2$};
\end{tikzpicture}
\end{gathered}
\quad \to \quad
\begin{gathered}
\begin{tikzpicture}[xscale=-0.5, yscale=0.5]
\draw [line,->-=1] (0,0) -- (0,-.282) node {$\times$} -- (0,-1.41) node [below] {$\cL_5$};
\draw [line,->-=.56] (0,0) -- (.5,.5) node [below left=-3pt] {$\cL_9$} -- (.8,.8) node {$+$} -- (1,1);
\draw [line,->-=.56] (1,1) -- (1.5,1.5) node [below left=-3pt] {$\cL_8$} -- (1.8,1.8) node {$+$} -- (2,2);
\draw [line,->-=1] (2,2) -- (3,3) node [above left=-3pt] {$\cL_4$};
\draw [line,->-=1] (0,0) -- (-3,3) node [above right=-3pt] {$\cL_1$};
\draw [line,->-=1] (1,1) -- (-1,3) node [above right=-3pt] {$\cL_2$};
\draw [line,->-=1] (2,2) -- (1,3) node [above right=-3pt] {$\cL_3$};
\end{tikzpicture}
\end{gathered} \, .
\]
This consistency condition is the {\it pentagon identity}
\ie
\label{Pentagon}
(F^{\cL_6,\cL_3,\cL_4}_{\ocL_5})_{\cL_7,\cL_8} \, (F^{\cL_1,\cL_2,\cL_8}_{\ocL_5})_{\cL_6,\cL_9}
\ = \ \sum_\cL \, (F^{\cL_1,\cL_2,\cL_3}_{\cL_7})_{\cL_6, \cL} \, (F^{\cL_1,\cL,\cL_4}_{\ocL_5})_{\cL_7,\cL_9} \, (F^{\cL_2,\cL_3,\cL_4}_{\cL_9})_{\cL,\cL_8} \, .
\fe
A solution to the pentagon identity amounts to the construction of a pivotal fusion category.  If there are $n$ isomorphism classes of simple objects, then the pentagon identity is a set of $\cO(n^9)$ cubic polynomial equations for $\cO(n^6)$ variables, modulo $\cO(n^3)$ gauge freedom.  As $n$ grows, a generic system of this size quickly becomes unmanageable.

The {\it cyclic permutation map} is the isomorphism relating the three vector spaces
\[
V_{\cL_1, \cL_2, \cL_3} \, , \quad V_{\cL_2, \cL_3, \cL_1} \, , \quad V_{\cL_3, \cL_1, \cL_2} \, ,
\]
which pictorially corresponds to moving the $\times$ mark around.  It is the $F$-move with an external edge representing the unit object $\cI$:
\[
\begin{gathered}
\begin{tikzpicture}[scale=1]
\draw [line,-<-=.56] (-1,0) -- (-.8,0) node {$\times$} -- (-.5,0) node [above] {$\ocL_3$} -- (0,0);
\draw [line,->-=1] (-1,0) -- (-1.5,.87) node [above left=-3pt] {$\cL_1$};
\draw [line,->-=1] (-1,0) -- (-1.5,-.87) node [below left=-3pt] {$\cL_2$};
\draw [line,->-=1] (0,0) -- (.5,-.87) node [below right=-3pt] {$\cL_3$};
\draw [line,dashed] (0,0) 
-- (.5,.87) node [above right=-3pt] {$\cI$};
\end{tikzpicture}
\end{gathered}
\quad = \quad
(F^{\cL_1, \cL_2, \cL_3}_{\cI})_{\ocL_3, \ocL_1}
\begin{gathered}
\begin{tikzpicture}[scale=1]
\draw [line,-<-=.56] (0,-1) -- (0,-.8) node {$\times$} -- (0,-.5) node [right] {$\ocL_1$} -- (0,0);
\draw [line,->-=1] (0,-1) -- (.87,-1.5) node [below right=-3pt] {$\cL_3$};
\draw [line,->-=1] (0,-1) -- (-.87,-1.5) node [below left=-3pt] {$\cL_2$};
\draw [line,->-=1] (0,0) -- (-.87,.5) node [above left=-3pt] {$\cL_1$};
\draw [line,dashed] (0,0) 
-- (.87,.5) node [above right=-3pt] {$\cI$};
\end{tikzpicture}
\end{gathered} \, .
\]
The net effect is a counter-clockwise rotation of the $\times$ mark accompanied by a factor of $(F^{\cL_1, \cL_2, \cL_3}_{\cI})_{\ocL_3, \ocL_1}$.  Gauge freedom alone cannot guarantee that the $F$-symbols $(F^{\cL_1, \cL_2, \cL_3}_{\cI})_{\ocL_3, \ocL_1}$ take value one for all $\cL_1, \cL_2, \cL_3$.\footnote{See for instance Appendix A of \cite{Chang:2018iay}.
}
The temptation to ignore the ordering and marking at trivalent vertices motivates the following definition.

\begin{defi}[Cyclic-permutation invariance]
A pivotal fusion category is called {\it cyclic-permutation invariant} if the trivalent vertices are cyclic-permutation invariant, that is, for every triple $(\cL_1, \cL_2, \cL_3)$ of objects,\footnote{A more conventional string diagram is
\[
\begin{gathered}
\begin{tikzpicture}[scale=.75]
\draw [line,->-=1] (0,0) -- (0,.5) node {$\times$} -- (0,1) node [above] {$\cL_3$};
\draw [line,->-=1] (0,0) -- (-.87,-.5) node [below left] {$\cL_1$};
\draw [line,->-=1] (0,0) -- (.87,-.5) node [below right] {$\cL_2$};
\end{tikzpicture}
\end{gathered}
\ = \
\begin{gathered}
\begin{tikzpicture}[scale=.75]
\draw [line,->-=1] (0,0) -- (0,.5) node {$\times$} -- (0,1) -- (1.73,0) node [below right] {$\cL_2$};
\draw [line,->-=1] (0,0) -- (-.87,-.5) -- (-1.73,0) node [above left] {$\cL_3$};
\draw [line,->-=1] (0,0) -- (.87,-.5) node [below right] {$\cL_1$};
\draw [line,dashed] (-.87,-.5) -- (-.87,-1) node [below] {$\cI$};
\draw [line,dashed] (0,1) -- (0,1.5) node [above] {$\cI$};
\end{tikzpicture}
\end{gathered}
\ = \
\begin{gathered}
\begin{tikzpicture}[xscale=-.75,yscale=.75]
\draw [line,->-=1] (0,0) -- (0,.5) node {$\times$} -- (0,1) -- (1.73,0) node [below left] {$\cL_1$};
\draw [line,->-=1] (0,0) -- (-.87,-.5) -- (-1.73,0) node [above right] {$\cL_3$};
\draw [line,->-=1] (0,0) -- (.87,-.5) node [below left] {$\cL_2$};
\draw [line,dashed] (-.87,-.5) -- (-.87,-1) node [below] {$\cI$};
\draw [line,dashed] (0,1) -- (0,1.5) node [above] {$\cI$};
\end{tikzpicture}
\end{gathered}
\]
involving evaluation, coevaluation, unitor and counitor.
}
\[
\begin{gathered}
\begin{tikzpicture}[scale=1]
\draw [line,->-=1] (0,0) -- (0,.5) node {$\times$} -- (0,1) node [above] {$\cL_3$};
\draw [line,->-=1] (0,0) -- (-.87,-.5) node [below left] {$\cL_1$};
\draw [line,->-=1] (0,0) -- (.87,-.5) node [below right] {$\cL_2$};
\end{tikzpicture}
\end{gathered}
\ = \
\begin{gathered}
\begin{tikzpicture}[scale=1]
\draw [line,->-=1] (0,0) -- (0,1) node [above] {$\cL_3$};
\draw [line,->-=1] (0,0) -- (-.87,-.5) node [below left] {$\cL_1$};
\draw [line,->-=1] (0,0) -- (.43,-.25) node {\rotatebox[origin=c]{60}{$\times$}} -- (.87,-.5) node [below right] {$\cL_2$};
\end{tikzpicture}
\end{gathered}
\ = \
\begin{gathered}
\begin{tikzpicture}[scale=1]
\draw [line,->-=1] (0,0) -- (0,1) node [above] {$\cL_3$};
\draw [line,->-=1] (0,0) -- (-.43,-.25) node {\rotatebox[origin=c]{30}{$\times$}} -- (-.87,-.5) node [below left] {$\cL_1$};
\draw [line,->-=1] (0,0) -- (.87,-.5) node [below right] {$\cL_2$};
\end{tikzpicture}
\end{gathered} \, .
\]
\end{defi}

In a cyclic-permutation invariant fusion category $\cC$, it is clear by a $\pi$-rotation that the $F$-symbols enjoy an order-two invariance
\[
(F^{\cL_1, \cL_2, \cL_3}_{\ocL_4})_{\cL_5, \cL_6}  = (F^{\cL_3, \cL_4, \cL_1}_{\ocL_2})_{\ocL_5, \ocL_6} \, .
\]
By relating the $F$-symbols to tetrahedra (as shown in Appendix~\ref{Sec:FTetra}), 
\ie
\label{FT}
&
\hspace{-.1in}
\begin{gathered}
\begin{tikzpicture}[scale=2]
\draw [line,->-=.18,-<-=.51,-<-=.83] (-1,0) -- (-.25,.43) node [above left=-2pt] {$\cL_1$} -- (.5,.87) -- (.5,0) node [right] {$\cL_6$} -- (.5,-.87) -- (-.25,-.43) node [below left=-2pt] {$\cL_2$} -- (-1,0);
\draw [line,->-=.5] (0,0) -- (-.33,0) node [above] {$\cL_5$} -- (-1,0);
\draw [line,->-=.5] (0,0) -- (.17,.29) node [above left=-2pt] {$\cL_4$} -- (.5,.87);
\draw [line,->-=.5] (0,0) -- (.17,-.29) node [below left=-2pt] {$\cL_3$} -- (.5,-.87);
\end{tikzpicture}
\end{gathered}
\ = \ (F^{\cL_1, \cL_2, \cL_3}_{\ocL_4})_{\cL_5, \cL_6}
\ 
{
\begin{gathered}
\begin{tikzpicture}[scale=1]
\begin{scope}
\node at (-1.4,0) {$\cL_1$};
\node at (1.4,0) {$\cL_6$};
\draw [line,->-=.01,-<-=.51] (0,0) circle (1);
\draw [line,->-=.53] (0,-1) -- (0,0) node [right] {$\cL_4$} -- (0,1);
\end{scope}
\begin{scope}[shift = {(3.75,0)}]
\node at (-1.4,0) {$\cL_6$};
\node at (1.4,0) {$\cL_2$};
\draw [line,->-=.01,->-=.51] (0,0) circle (1);
\draw [line,->-=.53] (0,-1) -- (0,0) node [right] {$\cL_3$} -- (0,1);
\end{scope}
\end{tikzpicture} 
\end{gathered} 
\over
\begin{gathered}
\begin{tikzpicture}[scale=1]
\node at (0,1.1) {};
\node at (0,0) {$\cL_6$};
\draw [line,-<-=.25] (0,0) circle (1);
\end{tikzpicture} 
\end{gathered}
}
\, ,
\\
&
\hspace{-.1in}
\begin{gathered}
\begin{tikzpicture}[xscale=-2, yscale=2]
\draw [line,->-=.18,-<-=.51,-<-=.83] (-1,0) -- (-.25,.43) node [above right=-2pt] {$\cL_4$} -- (.5,.87) -- (.5,0) node [left] {$\cL_6$} -- (.5,-.87) -- (-.25,-.43) node [below right=-2pt] {$\cL_3$} -- (-1,0);
\draw [line,-<-=.5] (0,0) -- (-.33,0) node [above] {$\cL_5$} -- (-1,0);
\draw [line,->-=.5] (0,0) -- (.17,.29) node [above right=-2pt] {$\cL_1$} -- (.5,.87);
\draw [line,->-=.5] (0,0) -- (.17,-.29) node [below right=-2pt] {$\cL_2$} -- (.5,-.87);
\end{tikzpicture}
\end{gathered}
\ = \ (F^{\cL_1, \cL_2, \cL_3}_{\ocL_4})_{\cL_5, \cL_6}
\ 
{
\begin{gathered}
\begin{tikzpicture}[scale=1]
\begin{scope}
\node at (-1.4,0) {$\cL_4$};
\node at (1.4,0) {$\cL_6$};
\draw [line,-<-=.01,->-=.51] (0,0) circle (1);
\draw [line,-<-=.53] (0,-1) -- (0,0) node [right] {$\cL_1$} -- (0,1);
\end{scope}
\begin{scope}[shift = {(3.75,0)}]
\node at (-1.4,0) {$\cL_6$};
\node at (1.4,0) {$\cL_3$};
\draw [line,-<-=.01,-<-=.51] (0,0) circle (1);
\draw [line,-<-=.53] (0,-1) -- (0,0) node [right] {$\cL_2$} -- (0,1);
\end{scope}
\end{tikzpicture} 
\end{gathered} 
\over
\begin{gathered}
\begin{tikzpicture}[scale=1]
\node at (0,1.1) {};
\node at (0,0) {$\cL_6$};
\draw [line,-<-=.25] (0,0) circle (1);
\end{tikzpicture} 
\end{gathered}
}
\, ,
\fe
additional relations can be manifested.  Each tetrahedron enjoys an $S_3$ symmetry: it is invariant under the $\bZ_3$ rotations and complex conjugate under a reflection.  Combined with the aforementioned $\pi$-rotation invariance generates an $S_4$ worth of relations for the $F$-symbols.  However, these relations are generally nonlinear due to the factors of graphs appearing on the right of \eqref{FT}.

\section{Transparent fusion categories}
\label{Sec:Transparent}

\begin{defi}[Transparency]
\label{Transparent}
A pivotal fusion category $\cC$ is called transparent if the associator involving any invertible object is the identity map.  In terms of string diagrams, $\cC$ is transparent if for every triple $(\cL_1, \cL_2, \cL_3)$ of objects in $\cC$ and for every invertible object $\eta$,
\[
\begin{gathered}
\begin{tikzpicture}[scale=1]
\draw [line,-<-=.56] (-1,0) 
-- (-.8,0) node {$\times$} 
-- (-.5,0) node [above] {$\overline{\eta\cL_3}$} -- (0,0);
\draw [line,->-=1] (-1,0) -- (-1.5,.87) node [above left=-3pt] {$\cL_1$};
\draw [line,->-=1] (-1,0) -- (-1.5,-.87) node [below left=-3pt] {$\cL_2$};
\draw [line,->-=1] (0,0) -- (.5,-.87) node [below right=-3pt] {$\cL_3$};
\draw [line,->-=1,dashed] (0,0) 
-- (.1,.175) node {\rotatebox[origin=c]{60}{$\times$}} 
-- (.5,.87) node [above right=-3pt] {$\eta$};
\end{tikzpicture}
\end{gathered}
\quad = \quad
\begin{gathered}
\begin{tikzpicture}[scale=1]
\draw [line,-<-=.56] (0,-1) 
-- (0,-.8) node {$\times$} 
-- (0,-.5) node [right] {$\overline{\cL_1\eta}$} -- (0,0);
\draw [line,->-=1] (0,-1) -- (.87,-1.5) node [below right=-3pt] {$\cL_3$};
\draw [line,->-=1] (0,-1) -- (-.87,-1.5) node [below left=-3pt] {$\cL_2$};
\draw [line,->-=1] (0,0) -- (-.87,.5) node [above left=-3pt] {$\cL_1$};
\draw [line,->-=1,dashed] (0,0) 
-- (.175,.1) node {\rotatebox[origin=c]{30}{$\times$}} 
-- (.87,.5) node [above right=-3pt] {$\eta$};
\end{tikzpicture}
\end{gathered}
\]
and likewise for $\eta$ on any of the three other external edges.
\end{defi}

Since the unit object is invertible, a transparent fusion category is automatically cyclic-permutation invariant.  Hence, the marking $\times$ on the trivalent vertices representing the ordering or edges can be ignored.  

Transparency essentially means that invertible objects can be attached or detached ``freely'', changing the isomorphism classes of the other involved objects without generating extra $F$-symbols.  Appendix~\ref{Sec:Graph} illustrates some basic operations.  The following operation will be referred to as {\it symmetry nucleation}:  given a graph, nucleate an invertible loop on any face and merge it with the bordering edges.  For example, on any triangular face,
\[
\begin{gathered}
\begin{tikzpicture}[scale=1]
\draw [line] (-1,0) -- (-1.43,-.25);
\draw [line] (1,0) -- (1.43,-.25);
\draw [line] (0,1.73) -- (0,2.23);
\draw [line,-<-=.53] (-1,0) -- (0,0) node [below] {$\cL_2$} -- (1,0);
\draw [line,->-=.53] (-1,0) -- (-.5,.87) node [above left=-2pt] {$\cL_1$} -- (0,1.73);
\draw [line,-<-=.53] (1,0) -- (.5,.87) node [above right=-2pt] {$\cL_3$} -- (0,1.73);
\draw [line,dashed,-<-=.27] (0,.58) circle (.375);
\node at (0,.58) {$\eta$};
\end{tikzpicture}
\end{gathered}
\quad = \quad
\begin{gathered}
\begin{tikzpicture}[scale=1]
\draw [line] (-1,0) -- (-1.43,-.25);
\draw [line] (1,0) -- (1.43,-.25);
\draw [line] (0,1.73) -- (0,2.23);
\draw [line,-<-=.53] (-1,0) -- (0,0) node [below] {$\eta\cL_2$} -- (1,0);
\draw [line,->-=.53] (-1,0) -- (-.5,.87) node [above left=-2pt] {$\eta\cL_1$} -- (0,1.73);
\draw [line,-<-=.53] (1,0) -- (.5,.87) node [above right=-2pt] {$\eta\cL_3$} -- (0,1.73);
\draw [line,dashed,->-=.6] (.66,.58) -- (.33,0);
\draw [line,dashed,->-=.6] (-.33,0) -- (-.66,.58);
\draw [line,dashed,->-=.6] (-.33,1.07) -- (.33,1.07);
\node at (0,.58) {$\eta$};
\end{tikzpicture}
\end{gathered}
\quad = \quad
\begin{gathered}
\begin{tikzpicture}[scale=1]
\draw [line] (-1,0) -- (-1.43,-.25);
\draw [line] (1,0) -- (1.43,-.25);
\draw [line] (0,1.73) -- (0,2.23);
\draw [line,-<-=.53] (-1,0) -- (0,0) node [below] {$\eta\cL_2$} -- (1,0);
\draw [line,->-=.53] (-1,0) -- (-.5,.87) node [above left=-2pt] {$\eta\cL_1$} -- (0,1.73);
\draw [line,-<-=.53] (1,0) -- (.5,.87) node [above right=-2pt] {$\eta\cL_3$} -- (0,1.73);
\end{tikzpicture}
\end{gathered} \, .
\]
A slight variant of symmetry nucleation gives rise to invariance relations for $F$-symbols.  Consider the $F$-move equation and add an invertible object $\eta$ to an open face
\ien
\begin{gathered}
\begin{tikzpicture}[scale=1]
\node at (-1.8,0) {$\eta$};
\draw [line,-<-=.56] (-1,0) -- (-.5,0) node [above] {$\cL_5$} -- (0,0);
\draw [line,->-=1] (-1,0) -- (-1.5,.87) node [above left=-3pt] {$\cL_1$};
\draw [line,->-=1] (-1,0) -- (-1.5,-.87) node [below left=-3pt] {$\cL_2$};
\draw [line,->-=1] (0,0) -- (.5,-.87) node [below right=-3pt] {$\cL_3$};
\draw [line,->-=1] (0,0) -- (.5,.87) node [above right=-3pt] {$\cL_4$};
\draw [line,dashed,->-=.55] (-2.5,0) ++(30:1) arc (30:-30:1);
\end{tikzpicture}
\end{gathered}
\quad &= \quad
\sum_{\cL} (F^{\cL_1, \cL_2, \cL_3}_{\ocL_4})_{\cL_5, \cL_6}
\begin{gathered}
\begin{tikzpicture}[scale=1]
\node at (-.8,-.5) {$\eta$};
\draw [line,-<-=.56] (0,-1) -- (0,-.5) node [right] {$\cL_6$} -- (0,0);
\draw [line,->-=1] (0,-1) -- (.87,-1.5) node [below right=-3pt] {$\cL_3$};
\draw [line,->-=1] (0,-1) -- (-.87,-1.5) node [below left=-3pt] {$\cL_2$};
\draw [line,->-=1] (0,0) -- (-.87,.5) node [above left=-3pt] {$\cL_1$};
\draw [line,->-=1] (0,0) -- (.87,.5) node [above right=-3pt] {$\cL_4$};
\draw [line,dashed,->-=.55] (-1.5,-.5) ++(30:1) arc (30:-30:1);
\end{tikzpicture}
\end{gathered} \, ,
\fen
which by transparency is equivalent to 
\ien
\begin{gathered}
\begin{tikzpicture}[scale=1]
\draw [line,-<-=.56] (-1,0) -- (-.5,0) node [above] {$\cL_5$} -- (0,0);
\draw [line,->-=1] (-1,0) -- (-1.5,.87) node [above left=-3pt] {$\cL_1\overline\eta$};
\draw [line,->-=1] (-1,0) -- (-1.5,-.87) node [below left=-3pt] {$\eta\cL_2$};
\draw [line,->-=1] (0,0) -- (.5,-.87) node [below right=-3pt] {$\cL_3$};
\draw [line,->-=1] (0,0) -- (.5,.87) node [above right=-3pt] {$\cL_4$};
\end{tikzpicture}
\end{gathered}
\quad = \quad
\sum_{\cL} (F^{\cL_1, \cL_2, \cL_3}_{\ocL_4})_{\cL_5, \cL_6}
\begin{gathered}
\begin{tikzpicture}[scale=1]
\draw [line,-<-=.56] (0,-1) -- (0,-.5) node [right] {$\eta\cL_6$} -- (0,0);
\draw [line,->-=1] (0,-1) -- (.87,-1.5) node [below right=-3pt] {$\cL_3$};
\draw [line,->-=1] (0,-1) -- (-.87,-1.5) node [below left=-3pt] {$\eta\cL_2$};
\draw [line,->-=1] (0,0) -- (-.87,.5) node [above left=-3pt] {$\cL_1\overline\eta$};
\draw [line,->-=1] (0,0) -- (.87,.5) node [above right=-3pt] {$\cL_4$};
\end{tikzpicture}
\end{gathered} \, .
\fen
The result is an invariance relation
\[
(F^{\cL_1, \cL_2, \cL_3}_{\ocL_4})_{\cL_5, \cL_6} \ = \ (F^{\cL_1\overline\eta, \eta\cL_2, \cL_3}_{\ocL_4})_{\cL_5, \eta\cL_6} \, .
\]
Similar operations on the other three faces give
\ien
(F^{\cL_1, \cL_2, \cL_3}_{\ocL_4})_{\cL_5, \cL_6} \ &= \ (F^{\cL_1, \cL_2\overline\eta, \eta\cL_3}_{\ocL_4})_{\cL_5\overline\eta, \cL_6} \
= \ (F^{\cL_1, \cL_2, \cL_3\overline\eta}_{\overline{\eta\cL_4}})_{\cL_5, \cL_6\overline\eta} \
= \ (F^{\eta\cL_1, \cL_2, \cL_3}_{\overline{\cL_4\overline\eta}})_{\eta\cL_5, \cL_6} \, .
\fen

Further useful relations between graphs and $F$-symbols can be derived as follows.  Let $(\cL_1, \cL_2, \cL_3)$ be any triple of simple objects in $\cC$, and $\eta$ any invertible object.  Consider 
\ie
\label{Mouse}
\begin{gathered}
\begin{tikzpicture}[scale=1]
\node at (0,2.4) {$\eta$};
\draw [line,dashed,->-=.51] (0,1) ++(-30:1) arc (-30:210:1);
\node at (-1.4,0) {$\cL_1$};
\node at (-.5,1.25) {$\cL_1\eta$};
\node at (1.4,0) {$\cL_3$};
\node at (.5,1.25) {$\cL_3\eta$};
\draw [line,-<-=.01,-<-=.18,-<-=.34,-<-=.51] (0,0) circle (1);
\draw [line,->-=.53] (0,-1) -- (0,0) node [right] {$\cL_2$} -- (0,1);
\end{tikzpicture}
\end{gathered}
\quad = \quad
\begin{gathered}
\begin{tikzpicture}[scale=1]
\node at (-1.4,0) {$\cL_1$};
\node at (1.4,0) {$\cL_3$};
\draw [line,-<-=.01,-<-=.51] (0,0) circle (1);
\draw [line,->-=.53] (0,-1) -- (0,0) node [right] {$\cL_2$} -- (0,1);
\end{tikzpicture}
\end{gathered}
 \, ,
\fe
and perform an $F$-move on $\cL_2$ to obtain
\ien
\begin{gathered}
\begin{tikzpicture}[scale=1]
\node at (0,2.4) {$\eta$};
\draw [line,dashed,->-=.51] (0,1) ++(-30:1) arc (-30:210:1);
\node at (-1.4,0) {$\cL_1$};
\node at (-.5,1.25) {$\cL_1\eta$};
\node at (1.4,0) {$\cL_3$};
\node at (.5,1.25) {$\cL_3\eta$};
\draw [line,-<-=.01,-<-=.18,-<-=.34,-<-=.51] (0,0) circle (1);
\draw [line,->-=.53] (0,-1) -- (0,0) node [right] {$\cL_2$} -- (0,1);
\end{tikzpicture}
\end{gathered}
\quad = \quad
(F^{\cL_1, \ocL_3, \cL_3\eta}
_{\cL_1\eta})_{\ocL_2, \eta}
\quad
\begin{gathered}
\begin{tikzpicture}[scale=1]
\begin{scope}
\node at (0,0) {$\cL_1$};
\draw [line,-<-=.25] (0,0) circle (1);
\end{scope}
\begin{scope}[shift = {(2.75,0)}]
\node at (0,0) {$\cL_3$};
\draw [line,-<-=.25] (0,0) circle (1);
\end{scope}
\end{tikzpicture} 
\end{gathered} \, .
\fen
Thus
\ie
\label{Graph1}
(F^{\cL_1, \ocL_3, \cL_3\eta}
_{\cL_1\eta})_{\ocL_2, \eta} \quad = \quad 
\frac{
\begin{gathered}
\begin{tikzpicture}[scale=1]
\node at (-1.4,0) {$\cL_1$};
\node at (1.4,0) {$\cL_3$};
\draw [line,-<-=.01,-<-=.51] (0,0) circle (1);
\draw [line,->-=.53] (0,-1) -- (0,0) node [right] {$\cL_2$} -- (0,1);
\end{tikzpicture}
\end{gathered}
}{
\begin{gathered}
\begin{tikzpicture}[scale=1]
\node at (0,1.1) {};
\begin{scope}
\node at (0,0) {$\cL_1$};
\draw [line,-<-=.25] (0,0) circle (1);
\end{scope}
\begin{scope}[shift = {(2.75,0)}]
\node at (0,0) {$\cL_3$};
\draw [line,-<-=.25] (0,0) circle (1);
\end{scope}
\end{tikzpicture} 
\end{gathered}
} \, .
\fe
The special case of $\cL_2 = \theta$ invertible, and $\cL_1 = \cL, ~ \cL_3 = \theta\cL$ gives
\ie
\label{Graph2}
{(F^{\cL, \overline{\theta\cL}, \theta\cL\eta}
_{\cL\eta})_{\overline\theta, \eta}}^{-1} \quad = \quad 
\begin{gathered}
\begin{tikzpicture}[scale=1]
\node at (0,1.1) {};
\begin{scope}
\node at (0,0) {$\cL$};
\draw [line,-<-=.25] (0,0) circle (1);
\end{scope}
\end{tikzpicture} 
\end{gathered}
\, .
\fe
Consider again the original diagram \eqref{Mouse}.  Perform an $F$-move on $\eta$, and then another $F$-move on a unit object connecting the two $\cL_2$ edges to obtain
\ien
&
\begin{gathered}
\begin{tikzpicture}[scale=1]
\node at (0,2.4) {$\eta$};
\draw [line,dashed,->-=.51] (0,1) ++(-30:1) arc (-30:210:1);
\node at (-1.4,0) {$\cL_1$};
\node at (-.5,1.25) {$\cL_1\eta$};
\node at (1.4,0) {$\cL_3$};
\node at (.5,1.25) {$\cL_3\eta$};
\draw [line,-<-=.01,-<-=.18,-<-=.34,-<-=.51] (0,0) circle (1);
\draw [line,->-=.53] (0,-1) -- (0,0) node [right] {$\cL_2$} -- (0,1);
\end{tikzpicture}
\end{gathered}
\quad = \quad
(F^{\ocL_1, \cL_1 \eta, \overline{\cL_3\eta}}
_{\ocL_3})_{\eta, \ocL_2}
\quad
\begin{gathered}
\begin{tikzpicture}[scale=1]
\draw [line,->-=.01,->-=.51] (0,0) circle (.5);
\node at (-1.9,0) {$\cL_1$};
\node at (-1,0) {$\cL_1\eta$};
\node at (1.9,0) {$\cL_3$};
\node at (1,0) {$\cL_3\eta$};
\draw [line,-<-=.01,-<-=.18,-<-=.34,-<-=.51] (0,0) circle (1.5);
\draw [line,->-=.55] (0,-1.5) -- (0,-1) node [right] {$\cL_2$} -- (0,-.5);
\draw [line,-<-=.55] (0,1.5) -- (0,1) node [right] {$\cL_2$} -- (0,.5);
\end{tikzpicture}
\end{gathered}
\\
\quad &= \quad
(F^{\ocL_1, \cL_1\eta, \overline{\cL_3\eta}}
_{\ocL_3})_{\eta, \ocL_2}
~~
(F^{\cL_2, \ocL_2, \cL_2}_{\ocL_2})_{\cI, \cI}
\quad
\begin{gathered}
\begin{tikzpicture}[scale=1]
\begin{scope}
\node at (-1.4,0) {$\cL_1$};
\node at (1.4,0) {$\cL_3$};
\draw [line,-<-=.01,-<-=.51] (0,0) circle (1);
\draw [line,->-=.53] (0,-1) -- (0,0) node [right] {$\cL_2$} -- (0,1);
\end{scope}
\begin{scope}[shift = {(3.75,0)}]
\node at (-1.5,0) {$\cL_1\eta$};
\node at (1.4,0) {$\cL_2$};
\draw [line,-<-=.01,->-=.51] (0,0) circle (1);
\draw [line,->-=.53] (0,-1) -- (0,0) node [right] {$\cL_3\eta$} -- (0,1);
\end{scope}
\end{tikzpicture} 
\end{gathered} \, .
\fen
From the above it is deduced that
\ie
\label{Graph3}
(F^{\ocL_1, \cL_1 \eta, \overline{\cL_3\eta}}
_{\ocL_3})_{\eta, \ocL_2}
\quad = \quad 
{
\begin{gathered}
\begin{tikzpicture}[scale=1]
\node at (0,0) {$\cL_2$};
\draw [line,-<-=.25] (0,0) circle (1);
\end{tikzpicture} 
\end{gathered}
\over
\begin{gathered}
\begin{tikzpicture}[scale=1]
\node at (0,1.1) {};
\node at (-1.5,0) {$\cL_1\eta$};
\node at (1.4,0) {$\cL_2$};
\draw [line,-<-=.01,->-=.51] (0,0) circle (1);
\draw [line,->-=.53] (0,-1) -- (0,0) node [right] {$\cL_3\eta$} -- (0,1);
\end{tikzpicture} 
\end{gathered} 
} \, .
\fe

\section{Transparent Haagerup-Izumi fusion categories}
\label{Sec:HI}

A Haagerup-Izumi fusion ring can be defined for every finite abelian group $G$.  A key feature is that it is quadratic \cite{thornton2012generalized,grossman2019infinite}: the fusion of a single non-invertible simple object with the invertible objects generate all the non-invertible simple objects.  In this section, a set of constraints are formulated for classifying transparent Haagerup-Izumi fusion categories with $G = \bZ_{2n+1}$.

The Haagerup-Izumi fusion ring with $G = \bZ_\nu$ has $\nu$ invertible objects
\[
\cI, \quad \A, \quad \A^2, \quad \cdots \quad \A^{\nu-1}
\]
and $\nu$ non-invertible simple objects
\[
\rho, \quad \A\rho, \quad \A^2\rho, \quad \cdots \quad \A^{\nu-1}\rho \, ,
\]
subject to the relations
\[
\A^\nu = 1 \, , \quad \A \rho = \rho \, \A^{\nu-1} \, , \quad \rho^2 = \cI + \sum_{k=0}^{\nu-1} \A^k \rho \, .
\]
When $\nu=1$, this is the Fibonacci ring, which is the Grothendieck ring of the Fibonacci category (even sectors of the $A_4$ subfactor) and Lee-Yang category.  When $\nu=2$, this is the Grothendieck ring of the $\cC(sl(2), 8)_{ad}$ fusion category (even sectors of the $A_7$ subfactor), which is premodular but not modular \cite{bruillard2016rank}.  When $\nu = 3$, this is the Grothendieck ring of the Haagerup $\cH_2$ and $\cH_3$ fusion categories \cite{asaeda1999exotic,Grossman_2012}.  For $\nu \ge 3$, the fusion ring is non-commutative.

Let $\cC$ be a transparent Haagerup-Izumi fusion category with $G = \bZ_{2n+1}$.  Define $\zeta$ and $\xi$ to be the graph values
\[
\zeta
\ \equiv \quad
\begin{gathered}
\begin{tikzpicture}[scale=1]
\node at (0,0) {$\rho$};
\draw [line] (0,0) circle (1);
\end{tikzpicture} 
\end{gathered} \, ,
\qquad
\xi
\ \equiv \quad
\begin{gathered}
\begin{tikzpicture}[scale=1]
\node at (-1.4,0) {$\rho$};
\node at (1.4,0) {$\rho$};
\draw [line] (0,0) circle (1);
\draw [line] (0,-1) -- (0,0) node [right] {$\rho$} -- (0,1);
\end{tikzpicture} 
\end{gathered} \, .
\]
On the left, symmetry nucleation implies that all non-invertible loops take value $\zeta$.  On the right, symmetry nucleation on the three faces implies that all such graphs with three non-invertible simple objects take the same value $\xi$.  In summary, for any triple $(\cL_1, \cL_2, \cL_3)$ of simple objects,
\[
\begin{gathered}
\begin{tikzpicture}[scale=1]
\node at (-1.4,0) {$\cL_1$};
\node at (1.4,0) {$\cL_3$};
\draw [line,->-=.01,-<-=.51] (0,0) circle (1);
\draw [line,->-=.53] (0,-1) -- (0,0) node [right] {$\cL_2$} -- (0,1);
\end{tikzpicture}
\end{gathered}
\quad = \quad 
\begin{cases}
1 & \text{all invertible} \, ,
\\
\zeta & \text{one invertible} \, ,
\\
\xi & \text{all non-invertible} \, .
\end{cases}
\]
By \eqref{Graph2}, for any pair $(\eta, \, \theta)$ of invertible objects,
\[
(F^{\cL, \overline{\theta\cL}, \theta\cL\eta}
_{\cL\eta})_{\overline\theta, \eta}
\ = \ 
\begin{cases}
1 & \cL ~~ \text{invertible} \, ,
\\
\zeta^{-1} & \cL ~~ \text{non-invertible} \, .
\end{cases}
\]

The $F$-symbols with a single internal invertible object can also be deduced.  For any triple $(\cL_1, \cL_2, \cL_3)$ of non-invertible simple objects, by \eqref{Graph1} and \eqref{Graph3},
\ie
& (F^{\cL_1, \cL_3, \eta\cL_3}_{\cL_1\overline{\eta}})_{\cL_2, \overline\eta} \ = \ \zeta^{-2} \, \xi \, ,
\qquad
(F^{\cL_1, \eta\cL_1, \eta\cL_3}_{\cL_3})_{\overline\eta, \cL_2} \ = \ \zeta \, \xi^{-1}
\, .
\fe
The possible values of $\zeta$ can be constrained as follows.  Consider two concentric $\rho$ loops and perform an $F$-move to obtain
\ien
\zeta^2 
\quad &= \quad
\begin{gathered}
\begin{tikzpicture}[scale=1]
\node at (0,0) {$\rho$};
\node at (-1.4,0) {$\rho$};
\draw [line] (0,0) circle (1);
\draw [line] (0,0) circle (.5);
\draw [line,dashed] (-.5,0) -- (-.75,0) node [above] {$\,\cI$} -- (-1,0);
\end{tikzpicture}
\end{gathered} 
\\
\quad &= \quad 
(F^{\rho, \rho, \rho}_{\rho})_{\cI, \cI} \quad
\begin{gathered}
\begin{tikzpicture}[scale=1]
\node at (-1.4,0) {$\cI$};
\node at (1.4,0) {$\rho$};
\draw [line] (0,0) ++(-90:1) arc (-90:90:1);
\draw [line,dashed] (0,0) ++(90:1) arc (90:270:1);
\draw [line] (0,-1) -- (0,0) node [right] {$\rho$} -- (0,1);
\end{tikzpicture}
\end{gathered}
\quad + \quad
\sum_{i = 0}^{2n}
\ (F^{\rho, \rho, \rho}_{\rho})_{\cI, \A^i\rho} \quad
\begin{gathered}
\begin{tikzpicture}[scale=1]
\node at (-1.4,0) {$\A^i\rho$};
\node at (1.4,0) {$\rho$};
\draw [line] (0,0) circle (1);
\draw [line] (0,0) circle (1);
\draw [line] (0,-1) -- (0,0) node [right] {$\rho$} -- (0,1);
\end{tikzpicture}
\end{gathered}
\\
\quad &= \quad 
1 + (2n+1) \, \zeta
\, .
\fen
Hence,
\[
\zeta \ = \ 
{2n+1 \pm \sqrt{(2n+1)^2+4} \over 2} \, .
\]
Finally, a gauge choice can be made such that
\[
\xi \ = \ \zeta^{\frac32} \, , 
\qquad
(F^{\cL_1, \cL_3, \eta\cL_3}_{\overline{\eta\cL_1}})_{\cL_2, \overline\eta}
\ = \ 
(F^{\cL_1, \eta\cL_1, \eta\cL_3}_{\ocL_3})_{\overline\eta, \cL_2}
\ = \
\zeta^{-\frac12} \, .
\]

\begin{defi}[Transparent constraints]
\label{T}
Let $I$ be the set of isomorphism classes of invertible objects and $N$ the set of isomorphism classes of non-invertible simple objects in the Haagerup-Izumi fusion ring with $G = \bZ_{2n+1}$.  The transparent constraints are the collection of constraints on the $F$-symbols
\ie
(F^{\eta, \cL_2, \cL_3}_{\cL_4})_{\cL_5, \cL_6} \ &= \ (F^{\cL_1, \eta, \cL_3}_{\cL_4})_{\cL_5, \cL_6}
\ = \ 
(F^{\cL_1, \cL_2, \eta}_{\cL_4})_{\cL_5, \cL_6}
\ = \ 
(F^{\cL_1, \cL_2, \cL_3}_{\overline\eta})_{\cL_5, \cL_6}
\ = \ 1 \, ,
\\
(F^{\eta\cL\theta, \overline{\theta}\cL, \cL}_{\eta\cL})_{\eta, \overline\theta}
\ &= \
\zeta^{-1} \, ,
\qquad
(F^{\cL_1, \cL_3, \eta\cL_3}_{\cL_1\overline\eta})_{\cL_2, \overline\eta}
\ = \ 
(F^{\cL_1, \eta\cL_1, \eta\cL_3}_{\cL_3})_{\overline\eta, \cL_2}
\ = \ 
\zeta^{-\frac12} \, ,
\\
(F^{\cL_1, \cL_2, \cL_3}_{\cL_4})_{\cL_5, \cL_6} \ &= \ (F^{\cL_1\overline\eta, \eta\cL_2, \cL_3}_{\cL_4})_{\cL_5, \eta\cL_6} \ 
= \ (F^{\cL_1, \cL_2\overline\eta, \eta\cL_3}_{\cL_4})_{\cL_5\overline\eta, \cL_6} \
\\
&= \ (F^{\cL_1, \cL_2, \cL_3\overline\eta}_{\cL_4\overline\eta})_{\cL_5, \cL_6\overline\eta} \
= \ (F^{\eta\cL_1, \cL_2, \cL_3}_{\eta\cL_4})_{\eta\cL_5, \cL_6} \, ,
\fe
for all $\eta, \theta \in I$ and $\cL, \cL_i \in N$.
\end{defi}

For the Haagerup-Izumi fusion ring with $G = \bZ_{2n+1}$, the number of independent $F$-symbols after imposing the transparent constraints is $(2n+1)^2+1$, significantly reduced from $\cO(n^6)$.  This number can be further reduced by exploiting tetrahedral invariance.  Since the factors in the relations \eqref{FT} between the tetrahedra and the $F$-symbols are universally equal to $\zeta^{-1} \xi^2$, the set of $F$-symbols with all objects non-invertible are {\it invariant} under the $A_4$ symmetry of the tetrahedron, and are related by complex conjugation under reflection if $\xi$ is chosen to be real.  To facilitate the computation, one may further assume reflection invariance and impose $S_4$ invariance on the $F$-symbols.\footnote{The usual notion of tetrahedral symmetry includes complex conjugation under reflections.  However, such relations complicate the present approach of computing a Groebner basis for the pentagon identity.  Hence, the term tetrahedral invariance in this paper refers to true equality without complex conjugation, and the notion is further subdivided into $A_4$ invariance (without reflections) and $S_4$ invariance (with reflections).
}

Table~\ref{Tab:Indep} lists the numbers of independent $F$-symbols after imposing the transparent constraints together with $A_4$ or $S_4$ tetrahedral invariance.  With $A_4$ invariance (necessary consequence of transparency), the pentagon identity under the transparent constraints can be practically solved up to $G = \bZ_9$ by computing a Groebner basis using MAGMA \cite{Magma}.  With $S_4$ invariance, it can be solved up to $G = \bZ_{15}$.  The next section presents the results of this classification.

\begin{table}[H]
\centering
\begin{tabular}{|c|c|c|}
\hline
$G$ & $A_4$ & $S_4$
\\\hline
$\bZ_3$ & 8 & 7
\\\hline
$\bZ_5$ & 22 & 16
\\\hline
$\bZ_7$ & 44 & 29
\\\hline
$\bZ_9$ & 74 & 46
\\\hline
$\bZ_{11}$ & 112 & 67
\\\hline
$\bZ_{13}$ & 158& 92
\\\hline
$\bZ_{15}$ & 212& 121
\\\hline
\end{tabular}
\caption{The numbers of independent $F$-symbols for the Haagerup-Izumi fusion rings after imposing the transparent constraints together with $A_4$ or $S_4$ tetrahedral invariance.}
\label{Tab:Indep}
\end{table}

\section{Classification of $F$-symbols}
\label{Sec:F}

\subsection{Main theorems}

\begin{thm}
\label{Main}
For the Haagerup-Izumi fusion rings with $G = \bZ_{2n+1}$, let
\[
\zeta_\pm \ \equiv \ {2n+1 \pm \sqrt{(2n+1)^2+4} \over 2} \, .
\]
Under the transparent constraints \eqref{T} and imposing $A_4$ tetrahedral invariance (necessary by transparency), the pentagon identity has the following solutions:
\begin{enumerate}[(a)]
\item There are two solutions for $G = \bZ_1$ corresponding to the Fibonacci and Lee-Yang categories.
\item There are eight solutions for $G = \bZ_3$.
\item There are sixteen solutions for $G = \bZ_5$.
\item There are twenty-four solutions for $G = \bZ_7$.
\item There are forty-eight solutions for $G = \bZ_9$.
\item For $G = \bZ_{2n+1}$ with $n = 1, 2, 3$, the solutions form four order-$2n$ orbits of the $\bZ_{2n}$ automorphism group.  Two orbits are unitary with $\zeta = \zeta_+$; the $F$-symbols are real in one of the two orbits, and complex in the other.  The remaining two orbits are the non-unitary Galois associates of the two unitary orbits, with $\zeta = \zeta_-$.  In particular, for $G = \bZ_3$, the unitary real orbit corresponds to the Haagerup $\cH_3$ fusion category, and the unitary complex orbit corresponds to the Haagerup $\cH_2$ fusion category, in the nomenclature of Grossman and Snyder \cite{Grossman_2012}.
\end{enumerate}
\end{thm}

\begin{thm}
\label{Main2}
For the Haagerup-Izumi fusion rings with $G = \bZ_{2n+1}$, let
\[
\zeta_\pm \ \equiv \ {2n+1 \pm \sqrt{(2n+1)^2+4} \over 2} \, .
\]  
Under the transparent constraints \eqref{T} and imposing $S_4$ tetrahedral invariance, the pentagon identity has the following solutions:
\begin{enumerate}[(a)]
\item There are two solutions for $G = \bZ_1$, corresponding to the Fibonacci and Lee-Yang categories.
\item There are four solutions for $G = \bZ_3$.
\item There are eight solutions for $G = \bZ_5$.
\item There are twelve solutions for $G = \bZ_7$.
\item There are twenty-four solutions for $G = \bZ_{13}$.
\item For $G = \bZ_{2n+1}$ with $n = 1, 2, 3, 6$, the solutions form two order-$2n$ orbits of the $\bZ_{2n}$ automorphism group.  One orbit is unitary with $\zeta = \zeta_+$, and the other orbit consists of the non-unitary Galois associates with $\zeta = \zeta_-$.  In particular, for $G = \bZ_3$, the unitary real orbit corresponds to the Haagerup $\cH_3$ fusion category in the nomenclature of Grossman and Snyder \cite{Grossman_2012}.
\item There are twenty-four solutions for $G = \bZ_9$, forming four order-six orbits of the $\bZ_6$ automorphism group.  Two orbits are unitary with $\zeta = \zeta_+$, and the other two orbits consist of the non-unitary Galois associates with $\zeta = \zeta_-$.
\item There are twenty-four solutions for $G = \bZ_{11}$, forming two order-two orbits and two order-ten orbits of the $\bZ_{10}$ automorphism group.  One order-two orbit and one order-ten orbit are unitary with $\zeta = \zeta_+$, and the other two orbits consist of the non-unitary Galois associates with $\zeta = \zeta_-$.
\item There are forty-eight solutions for $G = \bZ_{15}$, forming six order-eight orbits of the $\bZ_2 \times \bZ_4$ automorphism group.  Three orbits are unitary with $\zeta = \zeta_+$, and the other three orbits consist of the non-unitary Galois associates with $\zeta = \zeta_-$.
\item In the above, the $F$-symbols are real when $\zeta = \zeta_+$, and complex when $\zeta = \zeta_-$.  Solutions in a single orbit of the automorphism group have the same $(F^{\rho, \rho, \rho}_\rho)_{\rho, \rho}$, while different orbits have distinct $(F^{\rho, \rho, \rho}_\rho)_{\rho, \rho}$.  Since $(F^{\rho, \rho, \rho}_\rho)_{\rho, \rho}$ is gauge-invariant, solutions with distinct values correspond to inequivalent fusion categories.
\end{enumerate}
\end{thm}

\subsection{Explicit $F$-symbols for $G = \bZ_{2n+1}$ with $1 \le n \le 7$
}

Let $I$ be the set of invertible objects and $N$ the set of non-invertible simple objects of the Haagerup-Izumi fusion ring with $G = \bZ_{2n+1}$.  By \eqref{T}, the $F$-symbols involving at least one invertible object are given by
\[
(F^{\eta\cL\theta,\overline{\cL\theta}, \cL}_{\eta\cL})_{\eta, \overline\theta}
\ = \
\zeta^{-1} \, ,
\qquad
(F^{\cL_1, \cL_3, \eta\cL_3}_{\overline{\eta\cL_1}})_{\cL_2, \overline\eta}
\ = \ 
(F^{\cL_1, \eta\cL_1, \eta\cL_3}_{\ocL_3})_{\overline\eta, \cL_2}
\ = \ 
\zeta^{-\frac12} \, ,
\]
for all $\eta, \theta \in I$ and $\cL_i \in N$.  For the $F$-symbols with all simple objects being non-invertible, it suffices to specify the $(2n+1)^2$ components $(F^{\rho,\rho,\rho}_{*})_{\rho, *}$ with $*$ running over the non-invertible simple objects.  The rest are equal to one of the above by the $\bZ_{2n+1}^4$ invariance relations in \eqref{T}.  In fact, these invariance relations can be equivalently written as
\ien
(F_{\cL_4}^{\cL_1, \cL_2, \cL_3})_{\cL_5, \cL_6}
\ &= \ (F_{\eta\cL_4}^{\eta\cL_1, \eta\cL_2, \eta\cL_3})_{\eta\cL_5, \eta\cL_6}
\ = \ 
(F_{\cL_4}^{\eta\cL_1, \cL_2, \cL_3\eta})_{\cL_5, \cL_6}
\\
\ &= \ 
(F_{\cL_4\eta}^{\cL_1, \eta\cL_2, \cL_3})_{\cL_5, \cL_6}
\ = \ 
(F_{\cL_4}^{\cL_1, \cL_2, \cL_3})_{\eta\cL_5, \cL_6\eta} \, ,
\fen
for all $\eta, \theta \in I$ and $\cL_i \in N$.\footnote{Note that the equality of the first and the last terms implies that every $(F_{\cL_4}^{\cL_1, \cL_2, \cL_3})_{\cL_5, \cL_6}$ with $\cL_5, \, \cL_6 \in N$ is anti-circulant and therefore symmetric.  Together with the gauge choice in \eqref{T}, it follows that every $F$-symbol is symmetric in the appropriate basis.
}  

The explicit $F$-symbols for the Haagerup $\cH_2$ fusion category will first be presented, corresponding to two of the eight solutions in Theorem~\ref{Main}(b) that are unitary and complex (hence strictly $A_4$ tetrahedral invariance but not $S_4$).  Then, among the solutions classified by Theorem~\ref{Main2} (satisfying $S_4$ tetrahedral invariance), the explicit $F$-symbols for the unitary ones---the real ones with $\zeta = \zeta_+$---will now be presented.\footnote{The unitarity of the $F$-symbols for these solutions can be understood as follows.  First, by transparency, $F^{\cL_1,\cL_2, \cL_3}_{\cL_4}$ is nontrivial only if $\cL_1, \, \dotsc, \, \cL_4$ are all non-invertible and self-dual, so $F^{\ocL_1,\ocL_2, \ocL_3}_{\ocL_4} = F^{\cL_1,\cL_2, \cL_3}_{\cL_4}$.  Next, $S_4$ invariance and reality imply that $F^{\cL_2,\cL_3,\cL_4}_{\cL_1} = (F^{\ocL_1,\ocL_2, \ocL_3}_{\ocL_4})^t = (F^{\ocL_1,\ocL_2, \ocL_3}_{\ocL_4})^\dag$.  Hence unitarity $F^{\ocL_1,\ocL_2, \ocL_3}_{\ocL_4} (F^{\ocL_1,\ocL_2, \ocL_3}_{\ocL_4})^\dag = I$ becomes equivalent to the condition $F^{\cL_1,\cL_2, \cL_3}_{\cL_4} \, F^{\cL_2,\cL_3,\cL_4}_{\cL_1} = I$, which is the pentagon identity \eqref{Pentagon} with
$\cL_5 = \cI, \, \cL_7 = \cL_4 , \, \cL_9 = \cL_1$.
}
The $F$-symbols for the other fusion categories are available by request.  In the following, except for the particularly nice ones, the $F$-symbols will be presented as roots of polynomials, where $x_i$ denotes the $i$-th smallest real root of some polynomial in $x$, and likewise for other symbols $y, \, z, \, \dotsc$.  Note that all the presented polynomials only have simple roots.  The simpler polynomials are given in the main text, while the more complicated ones are given in Appendix~\ref{Sec:Poly}.

\subsubsection{Haagerup $\cH_2$ ($G = \bZ_3$)}

\paragraph{Theorem~\ref{Main}(b).}  Under the automorphism group $\text{Aut}(G) \cong \bZ_2$, there is exactly one unitary orbit with two complex solutions (strict $A_4$ tetrahedral invariance).  One solution is given by
\[
\begin{array}{c|ccc}
(F_{*}^{\rho,\rho,\rho})_{\rho, *} & \rho & \A\rho & \A^2\rho
\\\hline
\rho & x & z & z
\\
\A\rho & z & z & y_1
\\
\A^2\rho & z & y_2 & z
\end{array}
\]
where
\ien
& x = {7 - \sqrt{13} \over 6} \, , \quad y_{1, 2} = \frac{1}{12} \left(-1+\sqrt{13} \pm 3i\sqrt{2 \left(1+\sqrt{13}\right)}\right) \, , \quad z = {1 - \sqrt{13} \over 6} \, .
\fen
$\text{Aut}(G) \cong \bZ_2$ exchanges $y_1$ and $y_2$, giving the other solution in the orbit.

\subsubsection{Haagerup $\cH_3$ ($G = \bZ_3$)}

\paragraph{Theorem~\ref{Main2}(b).}  Under the automorphism group $\text{Aut}(G) \cong \bZ_2$, there is exactly one unitary orbit with two real solutions ($S_4$ tetrahedral invariance).  One solution is given by
\[
\begin{array}{c|ccc}
(F_{*}^{\rho,\rho,\rho})_{\rho, *} & \rho & \A\rho & \A^2\rho
\\\hline
\rho & x & y_1 & y_2
\\
\A\rho & y_1 & y_2 & z
\\
\A^2\rho & y_2 & z & y_1
\end{array}
\]
where
\ien
& x = {2 - \sqrt{13} \over 3} \, , \quad y_{1, 2} = \frac{1}{12} \left(5-\sqrt{13} \mp \sqrt{6 \left(1+\sqrt{13}\right)}\right) \, , \quad z = {1 + \sqrt{13} \over 6} \, .
\fen
$\text{Aut}(G) \cong \bZ_2$ exchanges $y_1$ and $y_2$, giving the other solution in the orbit.

\subsubsection{$G = \bZ_5$}

\paragraph{Theorem~\ref{Main2}(c).}  Under the automorphism group $\text{Aut}(G) \cong \bZ_4$, there is exactly one unitary orbit with four real solutions ($S_4$ tetrahedral invariance).  One solution is given by
\[
\begin{array}{c|ccccc}
(F_{*}^{\rho,\rho,\rho})_{\rho, *} & \rho & \A\rho & \A^2\rho & \A^3\rho & \A^4\rho
\\\hline
\rho & x & y_1 & y_3 & y_2 & y_4
\\
\A\rho & y_1 & y_4 & z_2 & z_4 & z_2
\\
\A^2\rho & y_3 & z_2 & y_2 & z_4 & z_4
\\
\A^3\rho & y_2 & z_4 & z_4 & y_3 & z_2
\\
\A^4\rho & y_4 & z_2 & z_4 & z_2 & y_1
\end{array}
\]
where
\[
x = {7 - \sqrt{29} \over 5} \, ,
\]
$y_i$ are the real roots of
\[
P^{\bZ_5}_y(y) = 625 y^8-1375 y^7+1275 y^6+245 y^5-654 y^4+152 y^3+75 y^2-29 y-1 \, ,
\]
and $z_i$ are the roots of
\[
P^{\bZ_5}_z(z) = 25 z^4-15 z^3-9 z^2+7 z-1 \, .
\]
$\text{Aut}(G) \cong \bZ_4$ permutes $y_i$ and exchanges $z_2$ and $z_4$ by
\ien
\tau_y=(1243) \, , \quad \tau_z=(24) \, ,
\fen
giving the other solutions in the orbit.

Note that the polynomial in $z$ factorizes over $\mathbb{Q}(\sqrt{29=5^2+4})$, and $z_2, \, z_4$ are the roots to one of the factors.  This pattern continues in the following solutions.  Namely, all polynomials factorize over $\mathbb{Q}(\sqrt{n^2+4})$, and the roots in a single orbit of $\text{Aut}(G)$ will always be roots of the same factor.

\subsubsection{$G = \bZ_7$}

\paragraph{Theorem~\ref{Main2}(d).}  Under the automorphism group $\text{Aut}(G) \cong \bZ_6$, there is exactly one unitary orbit with six solutions.  One solution is given by
\[
\begin{array}{c|ccccccc}
(F_{*}^{\rho,\rho,\rho})_{\rho, *} & \rho & \A\rho & \A^2\rho & \A^3\rho & \A^4\rho
& \A^5\rho& \A^6\rho
\\\hline
\rho & x & y_1 & y_2 & y_6 & y_4 & y_3 & y_5
\\
\A\rho & y_1 & y_5 & z_6 & w_2 & z_3 & w_1 & z_6
\\
\A^2\rho & y_2 & z_6 & y_3 & w_1 & z_4 & z_4 & w_2
\\
\A^3\rho & y_6 & w_2 & w_1 & y_4 & z_3 & z_4 & z_3
\\
\A^4\rho & y_4 & z_3 & z_4 & z_3 & y_6 & w_2 & w_1
\\
\A^5\rho & y_3 & w_1 & z_4 & z_4 & w_2 & y_2 & z_6
\\
\A^6\rho & y_5 & z_6 & w_2 & z_3 & w_1 & z_6 & y_1
\end{array}
\]
where
\[
x = {11 - 2\sqrt{53} \over 7} \, .
\]
$\text{Aut}(G) \cong \bZ_6 \cong \langle \sigma, \, \tau \mid \sigma^2 = \tau^3 = 1 \rangle$ 
permutes the roots by
\ien
& \sigma_y=(15)(23)(46) \, , \quad \sigma_z=\text{id} \, , \quad \sigma_w=(12) \, ,
\\
& \tau_y=(356)(142) \, , \quad \tau_z=(346) \, , \quad \tau_w=\text{id} \, ,
\fen
giving the other solutions in the orbit.

\subsubsection{$G = \bZ_9$}

\paragraph{Theorem~\ref{Main2}(g).}  Under the automorphism group $\text{Aut}(G) \cong \bZ_6$, there are two unitary orbits each with six solutions.  A solution in one orbit is given by
\[
\begin{array}{c|ccccccccc}
(F_{*}^{\rho,\rho,\rho})_{\rho, *} & \rho & \A\rho & \A^2\rho & \A^3\rho & \A^4\rho
& \A^5\rho& \A^6\rho& \A^7\rho& \A^8\rho
\\\hline
\rho     & x_1 & y_1 & y_{12} & \textcircled{r}_4 & y_6 & y_8 & \textcircled{r}_1 & y_7 & y_5
\\
\A\rho   & y_1 & y_5 & z_8 & w_{10} & w_2 & z_{11} & w_5 & w_7 & z_8
\\
\A^2\rho & y_{12} & z_8 & y_7 & w_7 & z_4 & w_9 & w_1 & z_4 & w_{10}
\\
\A^3\rho & \textcircled{r}_4 & w_{10} & w_7 & \textcircled{r}_1 & w_5 & w_1 & \textcircled{s}_4 & w_9 & w_2
\\
\A^4\rho & y_6 & w_2 & z_4 & w_5 & y_8 & z_{11} & w_9 & w_1 & z_{11}
\\
\A^5\rho & y_8 & z_{11} & w_9 & w_1 & z_{11} & y_6 & w_2 & z_4 & w_5
\\
\A^6\rho & \textcircled{r}_1 & w_5 & w_1 & \textcircled{s}_4 & w_9 & w_2 & \textcircled{r}_4 & w_{10} & w_7
\\
\A^7\rho & y_7 & w_7 & z_4 & w_9 & w_1 & z_4 & w_{10} & y_{12} & z_8
\\
\A^8\rho & y_5 & z_8 & w_{10} & w_2 & z_{11} & w_5 & w_7 & z_8 & y_1
\end{array}
\]
where
\[
x_{1, 2} ={35-4\sqrt{85} \mp \sqrt{517-56\sqrt{85}} \over 18}
\]
are the two negative roots of
\[
P_x^{\bZ_9}(x) = 81x^4-630x^3+899x^2+210x+9 \, .
\]
$\text{Aut}(G) \cong \bZ_6 \cong \langle \sigma, \, \tau \mid \sigma^2 = \tau^3 = 1 \rangle$ 
permutes the roots by
\ien
& \sigma_x=\text{id},\quad \sigma_y=(1\ 5)(2\ 4)(3\ 11)(6\ 8)(7\ 12)(9\ 10) \, , \quad \sigma_z=\text{id} \, ,
\\
& \sigma_w=(1\ 9)(2\ 5)(3\ 8)(4\ 12)(6\ 11)(7\ 10) \, , \quad \sigma_r=(1\ 4)(2\ 3) \, , \quad \sigma_s=\text{id} \, ,
\\
& \tau_x=\text{id},\quad \tau_y=(1\ 6\ 7)(2\ 3\ 9)(4\ 11\ 10)(5\ 8\ 12) \, , \quad \tau_z=(3\ 10\ 7)(4\ 8\ 11) \, ,
\\
& \tau_w=(1\ 7\ 2)(3\ 6\ 12)(4\ 8\ 11)(5\ 9\ 10) \, , \quad \tau_r=\text{id},\quad \tau_s=\text{id} \, ,
\fen
giving the other solutions in the orbit.
There is an additional $\bZ_2 \cong \langle \iota \mid \iota^2 = 1 \rangle$ action that acts by
\ien
& \iota_x=(1\ 2) \, , \quad \iota_y=(1\ 2)(3\ 6)(4\ 5)(7\ 9)(8\ 11)(10\ 12) \, , \quad \iota_z=(3\ 4)(7\ 11)(8\ 10) \, ,
\\ 
& \iota_w=(1\ 12)(2\ 6)(3\ 7)(4\ 9)(5\ 11)(8\ 10) \, , \quad \iota_r=(1\ 3)(2\ 4) \, , \quad \iota_s=(2\ 4) \, ,
\fen
and exchanges the two unitary orbits.

\subsubsection{$G = \bZ_{11}$}

\paragraph{Theorem~\ref{Main2}(h).} Under the automorphism group $\text{Aut}(G) \cong \bZ_{10}$, there is one unitary orbit with two solutions and one unitary orbit with ten solutions.  In the orbit with two solutions, one solution is given by
\[
\begin{array}{c|ccccccccccc}
(F_{*}^{\rho,\rho,\rho})_{\rho, *} & \rho & \A\rho & \A^2\rho & \A^3\rho & \A^4\rho
& \A^5\rho& \A^6\rho& \A^7\rho& \A^8\rho& \A^9\rho& \A^{10}\rho
\\\hline
\rho     & x & y_1 & y_2 & y_1 & y_1 & y_1 & y_2 & y_2 & y_2 & y_1 & y_2
\\
\A\rho   & y_1 & y_2 & z_2 & w_2 & w_2 & w_1 & z_2 & w_2 & w_1 & w_1 & z_2
\\
\A^2\rho & y_2 & z_2 & y_1 & w_1 & z_2 & w_2 & w_1 & w_2 & w_1 & z_2 & w_2
\\
\A^3\rho & y_1 & w_2 & w_1 & y_2 & w_1 & w_1 & z_2 & z_2 & z_2 & w_2 & w_2
\\
\A^4\rho & y_1 & w_2 & z_2 & w_1 & y_2 & w_2 & w_2 & z_2 & z_2 & w_1 & w_1
\\
\A^5\rho & y_1 & w_1 & w_2 & w_1 & w_2 & y_2 & z_2 & w_1 & z_2 & w_2 & z_1
\\
\A^6\rho & y_2 & z_2 & w_1 & z_2 & w_2 & z_2 & y_1 & w_1 & w_2 & w_1 & w_2
\\
\A^7\rho & y_2 & w_2 & w_2 & z_2 & z_2 & w_1 & w_1 & y_1 & w_2 & z_2 & w_1
\\
\A^8\rho & y_2 & w_1 & w_1 & z_2 & z_2 & z_2 & w_2 & w_2 & y_1 & w_2 & w_1
\\
\A^9\rho & y_1 & w_1 & z_2 & w_2 & w_1 & w_2 & w_1 & z_2 & w_2 & y_2 & z_2
\\
\A^{10}\rho & y_2 & z_2 & w_2 & w_2 & w_1 & z_1 & w_2 & w_1 & w_1 & z_2 & y_1
\end{array}
\]
where
\[
x=\frac{13-5\sqrt{5}}{11}
\]
is a root of the polynomial
\[
P_{2|x}^{\bZ_{11}}(x) = 11 x^2-26 x+4 \, .
\]
The $\bZ_2$ subgroup of $\text{Aut}(G) \cong \bZ_{10}$ exchanges $y_1$ with $y_2$ and $w_1$ with $w_2$.  In the order-ten orbit, one solution is given by
\[
\begin{array}{c|ccccccccccc}
(F_{*}^{\rho,\rho,\rho})_{\rho, *} & \rho & \A\rho & \A^2\rho & \A^3\rho & \A^4\rho
& \A^5\rho& \A^6\rho& \A^7\rho& \A^8\rho& \A^9\rho& \A^{10}\rho
\\\hline
\rho & x & y_1 & y_{10} & y_9 & y_2 & y_8 & y_3 & y_7 & y_4 & y_6 & y_5 \\
\A\rho &	y_1 & y_5 & z_6 & w_{10} & w_3 & w_9 & z_7 & w_1 & w_7 & w_4 & z_6 \\
\A^2\rho &	y_{10} & z_6 & y_6 & w_4 & z_3 & w_2 & w_6 & w_5 & w_8 & z_3 & w_{10} \\
\A^3\rho &	y_9 & w_{10} & w_4 & y_4 & w_7 & w_8 & z_4 & z_8 & z_4 & w_2 & w_3 \\
\A^4\rho &	y_2 & w_3 & z_3 & w_7 & y_7 & w_1 & w_5 & z_8 & z_8 & w_6 & w_9 \\
\A^5\rho &	y_8 & w_9 & w_2 & w_8 & w_1 & y_3 & z_7 & w_6 & z_4 & w_5 & z_7 \\
\A^6\rho &	y_3 & z_7 & w_6 & z_4 & w_5 & z_7 & y_8 & w_9 & w_2 & w_9 & w_1 \\
\A^7\rho &	y_7 & w_1 & w_5 & z_8 & z_8 & w_6 & w_9 & y_2 & w_3 & z_3 & w_7 \\
\A^8\rho &	y_4 & w_7 & w_8 & z_4 & z_8 & z_4 & w_2 & w_3 & y_9 & w_{10} & w_4 \\
\A^9\rho &	y_6 & w_4 & z_3 & w_2 & w_6 & w_5 & w_9 & z_3 & w_{10} & y_{10} & z_6 \\
\A^{10}\rho &	y_5 & z_6 & w_{10} & w_3 & w_9 & z_7 & w_1 & w_7 & w_4 & z_6 & y_1 \\
\end{array}
\]
where
\[
x = \frac{101 - 49 \sqrt5}{22}
\]
is a root of the polynomial
\[ 
P_{10|x}^{\bZ_{11}}(x) = 11 x^2-101 x-41 \, .
\]
$\text{Aut}(G) \cong \bZ_{10}\cong\langle \sigma, \, \tau \mid \sigma^2 = \tau^5 = 1 \rangle$ permutes the roots by
\ien
& \sigma_y=(1\ 5)(2\ 7)(3\ 8)(4\ 9)(6\ 10)\, , \quad \sigma_z=\text{id} \, , \quad \sigma_w=(1\ 9)(2\ 8)(3\ 7)(4\ 10)(5\ 6) \, ,
\fen
\ien
& \tau_y=(1\ 2\ 8\ 6\ 9)(3\ 10\ 4\ 5\ 7) \, , \quad \tau_z=(3\ 4\ 6\ 8\ 7) \, , \quad \tau_w=(1\ 5\ 2\ 10\ 3)(4\ 7\ 9\ 6\ 8) \, ,
\fen
giving the other solutions in the orbit.

\subsubsection{$G = \bZ_{13}$}

\paragraph{Theorem~\ref{Main2}(e).}  Under the automorphism group $\text{Aut}(G) \cong \bZ_{12}$, there is exactly one unitary orbit with twelve solutions.  One solution is given by
\[
\begin{array}{c|ccccccccccccc}
	(F_{*}^{\rho,\rho,\rho})_{\rho, *} & \rho & \A\rho & \A^2\rho & \A^3\rho & \A^4\rho
	& \A^5\rho& \A^6\rho& \A^7\rho& \A^8\rho& \A^9\rho& \A^{10}\rho& \A^{11}\rho& \A^{12}\rho
	\\\hline
\rho &	x & y_1 & y_9 & y_{12} & y_8 & y_4 & y_7 & y_3 & y_5 & y_2 & y_{10} & y_6 & y_{11} \\
\A\rho &	y_1 & y_{11} & z_6 & w_5 & s_3 & w_8 & w_4 & z_9 & w_{11} & w_9 & s_2 & w_2 & z_6 \\
\A^2\rho &	y_9 & z_6 & y_6 & w_2 & z_7 & w_{10} & w_{12} & s_1 & s_4 & w_7 & w_1 & z_7 & w_5 \\
\A^3\rho &	y_{12} & w_5 & w_2 & y_{10} & s_2 & w_1 & z_{10} & w_3 & z_8 & w_6 & z_{10} & w_{10} & s_3 \\
\A^4\rho &	y_8 & s_3 & z_7 & s_2 & y_2 & w_9 & w_7 & w_6 & z_4 & z_4 & w_3 & w_{12} & w_8 \\
\A^5\rho &	y_4 & w_8 & w_{10} & w_1 & w_9 & y_5 & w_{11} & s_4 & z_8 & z_4 & z_8 & s_1 & w_4 \\
\A^6\rho &	y_7 & w_4 & w_{12} & z_{10} & w_7 & w_{11} & y_3 & z_9 & s_1 & w_3 & w_6 & s_4 & z_9 \\
\A^7\rho &	y_3 & z_9 & s_1 & w_3 & w_6 & s_4 & z_9 & y_7 & w_4 & w_{12} & z_{10} & w_7 & w_{11} \\
\A^8\rho &	y_5 & w_{11} & s_4 & z_8 & z_4 & z_8 & s_1 & w_4 & y_4 & w_8 & w_{10} & w_1 & w_9 \\
\A^9\rho &	y_2 & w_9 & w_7 & w_6 & z_4 & z_4 & w_3 & w_{12} & w_8 & y_8 & s_3 & z_7 & s_2 \\
\A^{10}\rho &	y_{10} & s_2 & w_1 & z_{10} & w_3 & z_8 & w_6 & z_{10} & w_{10} & s_3 & y_{12} & w_5 & w_2 \\
\A^{11}\rho &	y_6 & w_2 & z_7 & w_{10} & w_{12} & s_1 & s_4 & w_7 & w_1 & z_7 & w_5 & y_9 & z_6 \\
\A^{12}\rho &	y_{11} & z_6 & w_5 & s_3 & w_8 & w_4 & z_9 & w_{11} & w_9 & s_2 & w_2 & z_6 & y_1 \\
\end{array}
\]
where
\[
x = \frac{107 - 8\sqrt{173}}{13}
\]
is a root of the polynomial
\[
P_x^{\bZ_{13}}(x) = 13 x^2-214 x+29 \, .
\]
$\text{Aut}(G) \cong \bZ_{12}\cong \langle \sigma, \,\tau \mid \sigma^4=\tau^3 = 1 \rangle$ permutes the roots in the following way
\ien
& \sigma_y=(1\ 4\ 11\ 5)(2\ 7\ 8\ 3)(6\ 12\ 9\ 10)\, , \quad \sigma_z=(4\ 9)(6\ 8)(7\ 10)\, , \\ &\quad \sigma_w=(1\ 5\ 10\ 2)(3\ 7\ 6\ 12)(4\ 9\ 11\ 8) \, ,\quad \sigma_s=(1\ 3\ 4\ 2)\, ,
\fen
\ien
& \tau_y=(1\ 2\ 12)(3\ 6\ 5)(4\ 7\ 9)(8\ 10\ 11) \, , \quad \tau_z=(4\ 10\ 6)(7\ 8\ 9) \, , \\ &\quad \tau_w=(1\ 4\ 7)(2\ 8\ 3)(5\ 9\ 6)(10\ 11\ 12) \, ,\quad \tau_s=\text{id}\, .
\fen

\subsubsection{$G = \bZ_{15}$}

\paragraph{Theorem~\ref{Main2}(i).}  Under the automorphism group $\text{Aut}(G) \cong \bZ_2 \times \bZ_4$, there are three unitary orbits each with eight solutions.  A solution in one orbit is given by\footnote{The polynomials are rather long and thus omitted in writing.
}
\[
\begin{array}{c|ccccccccccccccc}
(F_{*}^{\rho,\rho,\rho})_{\rho, *} & \rho & \A\rho & \A^2\rho & \A^3\rho & \A^4\rho
& \A^5\rho& \A^6\rho& \A^7\rho& \A^8\rho& \A^9\rho& \A^{10}\rho& \A^{11}\rho& \A^{12}\rho& \A^{13}\rho& \A^{14}\rho
\\\hline
\rho &x_2 & y_1 & y_9 & r_7 & y_2 & s_5 & r_1 & y_{23} & y_{16} & r_4 & s_4 & y_{17} & r_9 & y_{19} & y_{18} \\
\A \rho &y_1 & y_{18} & z_{14} & w_{14} & t_{10} & u_5 & v_{19} & w_{13} & z_7 & w_4 & v_{23} & u_{10} & t_5 & w_1 & z_{14} \\
\A^{2}\rho &y_9 & z_{14} & y_{19} & w_1 & z_{15} & v_6 & w_{10} & u_2 & t_7 & t_{12} & u_{12} & w_{20} & v_2 & z_{15} & w_{14} \\
\A^{3}\rho &r_7 & w_{14} & w_1 & r_9 & t_5 & v_2 & a_{11} & w_{18} & v_5 & a_4 & v_{10} & w_{19} & a_{11} & v_6 & t_{10} \\
\A^{4}\rho &y_2 & t_{10} & z_{15} & t_5 & y_{17} & u_{10} & w_{20} & w_{19} & z_4 & v_4 & v_{13} & z_4 & w_{18} & w_{10} & u_5 \\
\A^{5}\rho &s_5 & u_5 & v_6 & v_2 & u_{10} & s_4 & v_{23} & u_{12} & v_{10} & v_{13} & b_6 & v_4 & v_5 & u_2 & v_{19} \\
\A^{6}\rho &r_1 & v_{19} & w_{10} & a_{11} & w_{20} & v_{23} & r_4 & w_4 & t_{12} & a_4 & v_4 & v_{13} & a_4 & t_7 & w_{13} \\
\A^{7}\rho &y_{23} & w_{13} & u_2 & w_{18} & w_{19} & u_{12} & w_4 & y_{16} & z_7 & t_7 & v_5 & z_4 & v_{10} & t_{12} & z_7 \\
\A^{8}\rho &y_{16} & z_7 & t_7 & v_5 & z_4 & v_{10} & t_{12} & z_7 & y_{23} & w_{13} & u_2 & w_{18} & w_{19} & u_{12} & w_4 \\
\A^{9}\rho &r_4 & w_4 & t_{12} & a_4 & v_4 & v_{13} & a_4 & t_7 & w_{13} & r_1 & v_{19} & w_{10} & a_{11} & w_{20} & v_{23} \\
\A^{10}\rho &s_4 & v_{23} & u_{12} & v_{10} & v_{13} & z_{22} & v_4 & v_5 & u_2 & v_{19} & s_5 & u_5 & v_6 & v_2 & u_{10} \\
\A^{11}\rho &y_{17} & u_{10} & w_{20} & w_{19} & z_4 & v_4 & v_{13} & z_4 & w_{18} & w_{10} & u_5 & y_2 & t_{10} & z_{15} & t_5 \\
\A^{12}\rho &r_9 & t_5 & v_2 & a_{11} & w_{18} & v_5 & a_4 & v_{10} & w_{19} & a_{11} & v_6 & t_{10} & r_7 & w_{14} & w_1 \\
\A^{13}\rho &y_{19} & w_1 & z_{15} & v_6 & w_{10} & u_2 & t_7 & t_{12} & u_{12} & w_{20} & v_2 & z_{15} & w_{14} & y_9 & z_{14} \\
\A^{14}\rho &y_{18} & z_{14} & w_{14} & t_{10} & u_5 & v_{19} & w_{13} & z_7 & w_4 & v_{23} & u_{10} & t_5 & w_1 & z_{14} & y_1 \\
\end{array}
\]
$\text{Aut}(G) \cong \bZ_{2}\times \bZ_{4}\cong \langle \sigma, \,\tau\mid \sigma^2=\tau^4 = 1 \rangle$ permutes the roots in the following way
\ien
& \sigma_y=(1\ 18)(2\ 17)(9\ 19)(16\ 23)\, , \quad \sigma_r=(1\ 4)(7\ 9)\, , \quad \sigma_s=(4\ 5)\, , \quad \sigma_t=(5\ 10)(7\ 12)\, , \\&\sigma_r=(2\ 12)(5\ 10)\, , \quad \sigma_r=(1\ 4)(2\ 6)(4\ 13)(10\ 20)\, , \quad \sigma_w=(1\ 14)(4\ 13)(10\ 20)(18\ 19) \, ,\\&\sigma_z=\text{id}\, ,\quad \sigma_a=\text{id}\, ,\quad \sigma_b=\text{id}\, ,
\fen
\ien
& \tau_y=(1\ 19\ 2\ 23)(9\ 17\ 16\ 18)\, , \quad \tau_r=(1\ 7\ 4\ 9) \, , \quad \tau_s=\text{id} \, , \quad \tau_t=(5\ 7)(10\ 12) \, , \\&\tau_u=(2\ 10\ 12\ 5) \, , \quad \tau_v=(2\ 13\ 5\ 23)(4\ 10\ 19\ 6) \, , \quad \tau_w=(1\ 10\ 19\ 13)(4\ 14\ 20\ 18) \, ,\\&\tau_z=(4\ 7\ 14\ 15)\, ,\quad \tau_a=(4\ 11)\, ,\quad \tau_b=\text{id}\, ,
\fen
giving the other solutions in the orbit.  There is an additional $\bZ_3\cong \langle \iota|\iota^3=1\rangle$ action that acts by
\ien
& \iota_x=(1\ 2\ 3)\, , \quad \iota_y=(1\ 3\ 5)(2\ 10\ 14)(4\ 20\ 18)(6\ 11\ 17)(7\ 9\ 13)(8\ 12\ 19)(15\ 16\ 22)(21\ 24\ 23)\, , \\&\iota_r=(1\ 6\ 8)(2\ 7\ 5)(3\ 4\ 10)(9\ 12\ 11) \, , \quad \iota_s=(1\ 5\ 3)(2\ 6\ 4)\, , \\&\iota_t=(1\ 6\ 7)(2\ 5\ 8)(3\ 10\ 11)(4\ 12\ 9) \, , \quad \iota_u=(1\ 12\ 4)(2\ 6\ 7)(3\ 8\ 10)(5\ 11\ 9) \, , \\&\iota_v=(1\ 23\ 15)(2\ 3\ 12)(4\ 18\ 16)(5\ 7\ 22)(6\ 20\ 8)(9\ 14\ 13)(10\ 17\ 21)(11\ 24\ 19) \, , \\&\iota_w=(1\ 8\ 16)(2\ 18\ 11)(3\ 10\ 22)(4\ 9\ 6)(5\ 24\ 20)(7\ 21\ 13)(12\ 14\ 15)(17\ 19\ 23)\, ,\\&\iota_z=(4\ 5\ 11)(7\ 10\ 8)(14\ 19\ 17)(15\ 16\ 20)\, ,\\&\iota_a=(3\ 10\ 11)(4\ 8\ 6)\, ,\quad \iota_b=(3\ 6\ 5) \, ,
\fen
and cycles through the three distinct unitary orbits.  The polynomial for $x$ is given by 
\[
3375 x^{6}-116550 x^{5}+620280 x^{4}-926392 x^{3}+41520 x^{2}+88128 x-6912 \, .
\]

\section{Conclusions and outlook}
\label{Sec:Conclusions}

In this paper, the notion of a transparent fusion category was defined, and the $F$-symbols for transparent Haagerup-Izumi fusion categories with $G = \bZ_{2n+1}$ were constructively classified up to $G = \bZ_9$, and further up to $G = \bZ_{15}$ by additionally imposing $S_4$ tetrahedral invariance.  Various graph equivalences and $F$-symbol relations were derived from transparency, reducing the number of independent $F$-symbols from $\cO(n^6)$ to $\cO(n^2)$,ing the pentagon identity practically solvable.

In the Cuntz algebra approach to the construction of Haagerup-Izumi fusion categories, the polynomial equations are simpler to solve for Izumi systems \cite{izumi2001structure}, corresponding to $F$-symbols with strict $A_4$ tetrahedral invariance, than for Grossman-Snyder systems \cite{Grossman_2012}, corresponding to $F$-symbols with $S_4$ tetrahedral invariance.  For the former, the results of this paper (up to $G = \bZ_9$) can be directly compared with the solutions obtained by Evans and Gannon \cite{Evans:2010yr},and they are in complete agreement.  For the latter, whereas solutions to the Grossman-Snyder systems are available from Evans and Gannon \cite{Evans:2015zga} only up to $G = \bZ_5$, the present paper constructed categories up to $G = \bZ_{15}$.  To facilitate the comparison of the Cuntz algebra approach, the present authors solved the Grossman-Snyder system up to $G = \bZ_9$ by computing Groebner bases, and again found agreement.

Up to $G = \bZ_9$, the number of solutions under strict $A_4$ tetrahedral invariance (Izumi systems) and that under $S_4$ (Grossman-Snyder systems) are in agreement.  A possible explanation is that for any $G$, there is a one-to-one correspondence between the two, where a fusion category on one side is the bimodule category of another fusion category on the other side, with respect to an appropriate algebra object; in the physics language, they are related by gauging the $G = \bZ_{2n+1}$.  That this is true for the Haagerup \cite{haagerup1994principal} case $G = \bZ_3$ was shown by Grossman and Snyder \cite{Grossman_2012}.  If this explanation is generally valid, then for $G = \bZ_{15}$, the existence of three unitary orbits of the $\bZ_2 \times \bZ_4$ automorphism group according to Theorem~\ref{Main2} suggests that Evans and Gannon \cite{Evans:2010yr} missed two solutions in their numerical analysis.  However, the present authors have not been able to compute a Groebner basis for the $G = \bZ_{15}$ Izumi system to verify this speculation.

It would be interesting to construct transparent fusion categories for other fusion rings, especially quadratic (or generalized near-group) fusion rings where the fusion of the invertible objects with a single non-invertible object generates all the non-invertible objects \cite{thornton2012generalized,grossman2019infinite}.  Partial transparency may already be sufficient to reduce the pentagon identity to being practically solvable.  For instance, consider the following family of fusion rings:  introduce $\nu$ invertible objects
\[
\cI, \quad \A, \quad \A^2, \quad \cdots \quad \A^{\nu-1}
\]
together with $\nu+1$ non-invertible simple objects
\[
\rho, \quad \A\rho, \quad \A^2\rho, \quad \cdots \quad \A^{\nu-1}\rho, \quad \cN \, ,
\]
and define the fusion ring
\ien
& \A^\nu = 1 \, , \quad \A \rho = \rho \, \A^{\nu-1} \, , \quad \A \, \cN = \cN \, \A = \cN \, , 
\\
& \rho^2 = \cI + \cZ + \cN \, , \quad \cN^{\,2} = \cY + \cZ \, , \quad \rho \, \cN = \cN \rho = \cZ + \cN \, ,
\fen
where
\[
\cY \equiv \sum_{k=0}^{\nu-1} \A^k \, , \quad \cZ \equiv \sum_{k=0}^{\nu-1} \A^k \rho \, .
\]
When $\nu = 1$, this is the $R_\bC(\widehat{so}(3))_5$ fusion ring, which is known to admit three inequivalent fusion categories.  The generalization of $R_\bC(\widehat{so}(3))_5$ to the above family of fusion rings parallels the generalization of Fibonacci to Haagerup-Izumi.  Because the $\cN$ object is similar to the non-invertible object in the $G = \bZ_\nu$ Tambara-Yamagami categories, which are not transparent, it is unreasonable to expect that the above family extending the $R_\bC(\widehat{so}(3))_5$ fusion ring admits fully transparent fusion categories.  Nonetheless, partial transparency for $\rho, \, \A\rho, \, \dotsc, \, \A^{\nu-1} \rho$ may already be sufficient to render the pentagon identity solvable.

Explicit $F$-symbols have interesting applications.  For instance, three-manifold invariants can be defined by the $F$-symbols alone without the need of braiding \cite{Barrett:1993ab,gelfand1996invariants}.  In physics, one could study the gapped phase of (1+1)$d$ quantum field theories with Haagerup-Izumi symmetries, by constructing (1+1)$d$ topological quantum field theories with the same symmetries, as was done in \cite{Chang:2018iay} for fusion categories of smaller ranks.  The statistical models of \cite{Aasen:2020jwb} and the associated anyon chains can also be explicitly constructed from the unitary $F$-symbols.  In conformal field theory, the crossing symmetry of defect four-point functions may produce universal bounds on the spectra via the conformal bootstrap \cite{Simmons-Duffin:2016gjk}.

\section*{Acknowledgements}

The authors are grateful to Yuji Tachikawa for initiating their interest in the Haagerup-Izumi fusion categories, and to Matthew Titsworth for sharing his explicit $F$-symbols for the Haagerup $\mathcal{H}_3$ fusion category.  The authors also thank Chi-Ming Chang, Terry Gannon, Tobias Osborne, Yuji Tachikawa and Yifan Wang for helpful comments and suggestions on the draft.  This material is based upon work supported by the U.S. Department of Energy, Office of Science, Office of High Energy Physics, under Award Number DE-SC0011632.  YL is supported by the Sherman Fairchild Foundation.

\appendix

\section{$F$-symbols and tetrahedra}
\label{Sec:FTetra}
 
This appendix contains a derivation of the relation \eqref{FT} between $F$-symbols and tetrahedra.  On both sides of the $F$-move equation \eqref{F}
\[
\begin{gathered}
\begin{tikzpicture}[scale=1]
\draw [line,-<-=.56] (-1,0) -- (-.5,0) node [above] {$\cL_5$} -- (0,0);
\draw [line,->-=.56] (-1,0) -- (-1.5,.87) node [above left=-3pt] {$\cL_1$};
\draw [line,->-=.56] (-1,0) -- (-1.5,-.87) node [below left=-3pt] {$\cL_2$};
\draw [line,->-=.56] (0,0) -- (.5,-.87) node [below right=-3pt] {$\cL_3$};
\draw [line,->-=.56] (0,0) -- (.5,.87) node [above right=-3pt] {$\cL_4$};
\end{tikzpicture}
\end{gathered}
\quad = \quad 
\sum_{\cL} (F^{\cL_1, \cL_2, \cL_3}_{\ocL_4})_{\cL_5, \cL}
\begin{gathered}
\begin{tikzpicture}[scale=1]
\draw [line,-<-=.56] (0,-1) -- (0,-.5) node [right] {$\cL$} -- (0,0);
\draw [line,->-=.56] (0,-1) -- (.87,-1.5) node [below right=-3pt] {$\cL_3$};
\draw [line,->-=.56] (0,-1) -- (-.87,-1.5) node [below left=-3pt] {$\cL_2$};
\draw [line,->-=.56] (0,0) -- (-.87,.5) node [above left=-3pt] {$\cL_1$};
\draw [line,->-=.56] (0,0) -- (.87,.5) node [above right=-3pt] {$\cL_4$};
\end{tikzpicture}
\end{gathered} \, ,
\]
join 
\ie
\label{I}
\begin{gathered}
\begin{tikzpicture}[scale=1]
\draw [line,->-=.56] (0,-1) -- (0,-.5) node [right] {$\cL_6$} -- (0,0);
\draw [line,-<-=.56] (0,-1) -- (.87,-1.5) node [below right=-3pt] {$\cL_2$};
\draw [line,-<-=.56] (0,-1) -- (-.87,-1.5) node [below left=-3pt] {$\cL_3$};
\draw [line,-<-=.56] (0,0) -- (-.87,.5) node [above left=-3pt] {$\cL_4$};
\draw [line,-<-=.56] (0,0) -- (.87,.5) node [above right=-3pt] {$\cL_1$};
\end{tikzpicture}
\end{gathered}
\fe
from the right.  The resulting graph on the left side of the $F$-move equation can be adjusted into a tetrahedron
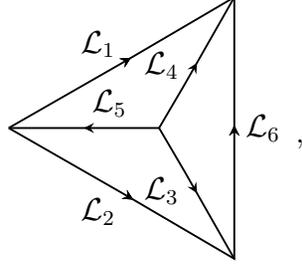
\[
\begin{gathered}
\begin{tikzpicture}[scale=2]
\draw [line,->-=.18,-<-=.51,-<-=.83] (-1,0) -- (-.25,.43) node [above left=-2pt] {$\cL_1$} -- (.5,.87) -- (.5,0) node [right] {$\cL_6$} -- (.5,-.87) -- (-.25,-.43) node [below left=-2pt] {$\cL_2$} -- (-1,0);
\draw [line,->-=.5] (0,0) -- (-.33,0) node [above] {$\cL_5$} -- (-1,0);
\draw [line,->-=.5] (0,0) -- (.17,.29) node [above left=-2pt] {$\cL_4$} -- (.5,.87);
\draw [line,->-=.5] (0,0) -- (.17,-.29) node [below left=-2pt] {$\cL_3$} -- (.5,-.87);
\end{tikzpicture}
\end{gathered} \, ,
\]
whereas the graph on the right side of the $F$-move equation can be adjusted into
\[
\begin{gathered}
\begin{tikzpicture}[scale=1]
\draw[line,->-=.76] (0,1) ellipse (1 and .5);
\draw[line,-<-=.26] (0,1) ellipse (1 and .5);
\draw[line,->-=.76] (0,-1) ellipse (1 and .5);
\draw[line,-<-=.26] (0,-1) ellipse (1 and .5);
\draw[line,->-=.5] (1,-1) -- (1,0) node [right] {$\cL_6$} -- (1,1);
\draw[line,-<-=.5] (-1,-1) -- (-1,0) node [left] {$\cL$}-- (-1,1);
\node at (0,.75) {$\cL_4$};
\node at (0,-1.25) {$\cL_2$};
\node at (0,-.25) {$\cL_3$};
\node at (0,1.75) {$\cL_1$};
\end{tikzpicture}
\end{gathered}
\ = \
\D_{\cL, \cL_6}
\times
\begin{gathered}
\begin{tikzpicture}[scale=1]
\draw[line,->-=.76] (0,1) ellipse (1 and .5);
\draw[line,-<-=.26] (0,1) ellipse (1 and .5);
\draw[line,->-=.76] (0,-1) ellipse (1 and .5);
\draw[line,-<-=.26] (0,-1) ellipse (1 and .5);
\draw[line,->-=.5] (1,-1) -- (1,0) node [right] {$\cL_6$} -- (1,1);
\draw[line,-<-=.5] (-1,-1) -- (-1,0) node [left] {$\cL_6$}-- (-1,1);
\node at (0,.75) {$\cL_4$};
\node at (0,-1.25) {$\cL_2$};
\node at (0,-.25) {$\cL_3$};
\node at (0,1.75) {$\cL_1$};
\end{tikzpicture}
\end{gathered} \, ,
\]
which vanishes if $\cL \neq \cL_6$ because the top and bottom loops can be shrunk but the vector space $V_{\cL, \cL_6}$ is empty.  Applying the $F$-move to a unit object connecting the two $\cL_6$ edges gives
\ien
\begin{gathered}
\begin{tikzpicture}[scale=1]
\draw[line,->-=.76] (0,1) ellipse (1 and .5);
\draw[line,-<-=.26] (0,1) ellipse (1 and .5);
\draw[line,->-=.76] (0,-1) ellipse (1 and .5);
\draw[line,-<-=.26] (0,-1) ellipse (1 and .5);
\draw[line,->-=.5] (1,-1) -- (1,0) node [right] {$\cL_6$} -- (1,1);
\draw[line,-<-=.5] (-1,-1) -- (-1,0) node [left] {$\cL_6$}-- (-1,1);
\node at (0,.75) {$\cL_4$};
\node at (0,-1.25) {$\cL_2$};
\node at (0,-.25) {$\cL_3$};
\node at (0,1.75) {$\cL_1$};
\end{tikzpicture}
\end{gathered}
\quad = \quad
(F^{\ocL_6, \cL_6, \ocL_6}_{\ocL_6})_{\cI, \cI}
\quad
\begin{gathered}
\begin{tikzpicture}[scale=1]
\begin{scope}
\node at (-1.4,0) {$\cL_1$};
\node at (1.4,0) {$\cL_6$};
\draw [line,->-=.01,-<-=.51] (0,0) circle (1);
\draw [line,->-=.53] (0,-1) -- (0,0) node [right] {$\cL_4$} -- (0,1);
\end{scope}
\begin{scope}[shift = {(3.75,0)}]
\node at (-1.4,0) {$\cL_6$};
\node at (1.4,0) {$\cL_2$};
\draw [line,->-=.01,->-=.51] (0,0) circle (1);
\draw [line,->-=.53] (0,-1) -- (0,0) node [right] {$\cL_3$} -- (0,1);
\end{scope}
\end{tikzpicture} 
\end{gathered} \, .
\fen
Again, no non-unit object $\cL$ can bridge the two $\Theta$ graphs on the right because the $\Theta$ graphs can be shrunk, but the vector space $V_{\cL,\cI}$ is empty if $\cL \neq \cI$.

Putting things together,
\ien
\hspace{-.1in}
\begin{gathered}
\begin{tikzpicture}[scale=2]
\draw [line,->-=.18,-<-=.51,-<-=.83] (-1,0) -- (-.25,.43) node [above left=-2pt] {$\cL_1$} -- (.5,.87) -- (.5,0) node [right] {$\cL_6$} -- (.5,-.87) -- (-.25,-.43) node [below left=-2pt] {$\cL_2$} -- (-1,0);
\draw [line,->-=.5] (0,0) -- (-.33,0) node [above] {$\cL_5$} -- (-1,0);
\draw [line,->-=.5] (0,0) -- (.17,.29) node [above left=-2pt] {$\cL_4$} -- (.5,.87);
\draw [line,->-=.5] (0,0) -- (.17,-.29) node [below left=-2pt] {$\cL_3$} -- (.5,-.87);
\end{tikzpicture}
\end{gathered}
\ = \ (F^{\cL_1, \cL_2, \cL_3}_{\ocL_4})_{\cL_5, \cL_6}
\ 
{
\begin{gathered}
\begin{tikzpicture}[scale=1]
\begin{scope}
\node at (-1.4,0) {$\cL_1$};
\node at (1.4,0) {$\cL_6$};
\draw [line,->-=.01,-<-=.51] (0,0) circle (1);
\draw [line,->-=.53] (0,-1) -- (0,0) node [right] {$\cL_4$} -- (0,1);
\end{scope}
\begin{scope}[shift = {(3.75,0)}]
\node at (-1.4,0) {$\cL_6$};
\node at (1.4,0) {$\cL_2$};
\draw [line,->-=.01,->-=.51] (0,0) circle (1);
\draw [line,->-=.53] (0,-1) -- (0,0) node [right] {$\cL_3$} -- (0,1);
\end{scope}
\end{tikzpicture} 
\end{gathered} 
\over
\begin{gathered}
\begin{tikzpicture}[scale=1]
\node at (0,1.1) {};
\node at (0,0) {$\cL_6$};
\draw [line,-<-=.25] (0,0) circle (1);
\end{tikzpicture} 
\end{gathered}
}
\, .
\fen
A similar derivation by joining \eqref{I} from the left with the $F$-move equation shows that
\ien
\hspace{-.1in}
\begin{gathered}
\begin{tikzpicture}[xscale=-2, yscale=2]
\draw [line,->-=.18,-<-=.51,-<-=.83] (-1,0) -- (-.25,.43) node [above right=-2pt] {$\cL_4$} -- (.5,.87) -- (.5,0) node [left] {$\cL_6$} -- (.5,-.87) -- (-.25,-.43) node [below right=-2pt] {$\cL_3$} -- (-1,0);
\draw [line,-<-=.5] (0,0) -- (-.33,0) node [above] {$\cL_5$} -- (-1,0);
\draw [line,->-=.5] (0,0) -- (.17,.29) node [above right=-2pt] {$\cL_1$} -- (.5,.87);
\draw [line,->-=.5] (0,0) -- (.17,-.29) node [below right=-2pt] {$\cL_2$} -- (.5,-.87);
\end{tikzpicture}
\end{gathered}
\ = \ (F^{\cL_1, \cL_2, \cL_3}_{\ocL_4})_{\cL_5, \cL_6}
\ 
{
\begin{gathered}
\begin{tikzpicture}[scale=1]
\begin{scope}
\node at (-1.4,0) {$\cL_4$};
\node at (1.4,0) {$\cL_6$};
\draw [line,-<-=.01,->-=.51] (0,0) circle (1);
\draw [line,-<-=.53] (0,-1) -- (0,0) node [right] {$\cL_1$} -- (0,1);
\end{scope}
\begin{scope}[shift = {(3.75,0)}]
\node at (-1.4,0) {$\cL_6$};
\node at (1.4,0) {$\cL_3$};
\draw [line,-<-=.01,-<-=.51] (0,0) circle (1);
\draw [line,-<-=.53] (0,-1) -- (0,0) node [right] {$\cL_2$} -- (0,1);
\end{scope}
\end{tikzpicture} 
\end{gathered} 
\over
\begin{gathered}
\begin{tikzpicture}[scale=1]
\node at (0,1.1) {};
\node at (0,0) {$\cL_6$};
\draw [line,-<-=.25] (0,0) circle (1);
\end{tikzpicture} 
\end{gathered}
}
\, .
\fen

\section{Transparent graph equivalences}
\label{Sec:Graph}

Let $\cC$ be a transparent fusion category, and $\eta$ an invertible object.  There are the following graph equivalences.

\begin{enumerate}
\item {\bf (Loop Value)} ~
Applying the $F$-move to an invertible $\eta$ loop gives
\ien
\begin{gathered}
\begin{tikzpicture}[scale=1]
\node at (-1.4,0) {$\eta$};
\node at (1.4,0) {$\eta$};
\draw [line,dashed,-<-=.01,-<-=.51] (0,0) circle (1);
\draw [line,dashed] (0,-1) -- (0,0) node [right] {$\cI$} -- (0,1);
\end{tikzpicture}
\end{gathered}
\quad = \quad
\begin{gathered}
\begin{tikzpicture}[scale=1]
\node at (-2,0) {$\eta$};
\node at (2,0) {$\eta$};
\draw [line,dashed,-<-=.25] (-2,0) circle (1);
\draw [line,dashed,-<-=.25] (2,0) circle (1);
\draw [line,dashed] (-1,0) -- (0,0) node [above] {$\cI$} -- (1,0);
\end{tikzpicture}
\end{gathered} \, .
\fen
Thus
\ie
\begin{gathered}
\begin{tikzpicture}[scale=1]
\node at (0,0) {$\eta$};
\draw [line,dashed,-<-=.25] (0,0) circle (1);
\end{tikzpicture}
\end{gathered} 
\quad = \ 1 \, ,
\fe
{\it i.e.} invertible loops have value 1.

\item {\bf (Attachment)} ~
An invertible object can be attached to a simple object $\cL$
\ien
\begin{gathered}
\begin{tikzpicture}[scale=1]
\node at (-.9,0) {$\eta$};
\draw [line,-<-=.53] (0,-1) -- (0,0) node [right] {$\cL$} -- (0,1);
\draw [line,dashed,->-=.55] (-1.5,0) ++(30:1) arc (30:-30:1);
\end{tikzpicture}
\end{gathered}
\quad = \quad
\begin{gathered}
\begin{tikzpicture}[scale=1]
\draw [line,-<-=.53] (0,-1) node [left] {$\cL$} -- (0,0) node [left] {$\eta\cL$} -- (0,1) node [left] {$\cL$};
\draw [line,dashed,->-=.5] (-1,.5) -- (0,.5);
\draw [line,dashed,-<-=.5] (-1,-.5) -- (0,-.5);
\end{tikzpicture}
\end{gathered}
\quad \, .
\fen

\item  {\bf (Detachment)} ~
An invertible object with two ends attached to a non-invertible simple object $\cL$ can be detached
\ien
\begin{gathered}
\begin{tikzpicture}[scale=1]
\node at (-.9,0) {$\eta$};
\draw [line,-<-=.53,-<-=.03] (0,-1) node [right] {$\cL$} -- (0,0) node [right] {$\eta\cL$} -- (0,1) node [right] {$\cL$};
\draw [line,dashed,->-=.55] (0,0) ++(270:.5) arc (270:90:.5);
\end{tikzpicture}
\end{gathered} 
\quad = \quad
\begin{gathered}
\begin{tikzpicture}[scale=1]
\node at (-2,0) {$\eta$};
\draw [line,-<-=.53] (0,-1) -- (0,0) node [right] {$\cL$} -- (0,1);
\draw [line,dashed,-<-=.5] (-2,0) circle (.5);
\draw [line,dashed] (-1.5,0) -- (-.75,0) node [above] {$\cI$} -- (0,0);
\end{tikzpicture}
\end{gathered}
\quad = \quad
\begin{gathered}
\begin{tikzpicture}[scale=1]
\draw [line,-<-=.53] (0,-1) -- (0,0) node [right] {$\cL$} -- (0,1);
\end{tikzpicture}
\end{gathered} \, .
\fen

\item {\bf (Swap)} ~
An invertible object attached to an edge can be swapped across a trivalent vertex
\ien
&
\begin{gathered}
\begin{tikzpicture}[scale=1]
\draw [line,-<-=.53,-<-=.03] (0,-1) node [left] {$\cL_1$} -- (0,0) node [left] {$\eta\cL_1$} -- (0,1) node [left] {$\cL_3$};
\draw [line,-<-=.53,-<-=.03] (1,.5) node [right] {$\cL_2$} -- (0,.5);
\draw [line,dashed,-<-=.5] (-1,-.5) node [left] {$\eta$} -- (0,-.5);
\end{tikzpicture}
\end{gathered}
\quad = \quad
\begin{gathered}
\begin{tikzpicture}[scale=1]
\draw [line,-<-=.53] (0,-1) node [left] {$\cL_1$} -- (0,0) node [left] {$\overline\eta\cL_3$} -- (0,1) node [left] {$\cL_3$};
\draw [line,-<-=.53] (1,-.5) node [right] {$\cL_2$} -- (0,-.5);
\draw [line,dashed,-<-=.5] (-1,.5) node [left] {$\eta$} -- (0,.5);
\end{tikzpicture}
\end{gathered} \, .
\fen

\item {\bf (Contraction)} ~
An invertible object bridged across a trivalent vertex can be contracted.  It can be regarded as a swap followed by a detachment
\ien
\begin{gathered}
\begin{tikzpicture}[scale=1]
\draw [line,->-=1] (0,0) -- (-1,0) node [left] {$\cL_1$};
\draw [line,->-=1] (0,0) node [below=5] {$\eta\cL_2~~~$} node [above=8] {$\eta\cL_3~~~$} -- (.5,-.87) node [below right] {$\cL_2$};
\draw [line,->-=1] (0,0) -- (.5,.87) node [above right] {$\cL_3$};
\draw [line,dashed,-<-=.53] (.33,-.58) -- (.33,0) node [right] {$\eta$} -- (.33,.58);
\end{tikzpicture}
\end{gathered}
\quad = \quad
\begin{gathered}
\begin{tikzpicture}[scale=1]
\draw [line,->-=1] (0,0) -- (-1,0) node [left] {$\cL$};
\draw [line,->-=1] (0,0) -- (.5,-.87) node [below right] {$\cL_2$};
\draw [line,->-=1] (0,0) -- (.5,.87) node [above right] {$\cL_3$};
\draw [line,dashed,-<-=.53] (.225,.39) ++(-120:.3) arc (-120:60:.3);
\node at (.75,0) {$\eta$};
\end{tikzpicture}
\end{gathered}
\quad &= \quad
\begin{gathered}
\begin{tikzpicture}[scale=1]
\draw [line,->-=1] (0,0) -- (-1,0) node [left] {$\cL_1$};
\draw [line,->-=1] (0,0) -- (.5,-.87) node [below right] {$\cL_2$};
\draw [line,->-=1] (0,0) -- (.5,.87) node [above right] {$\cL_3$};
\end{tikzpicture}
\end{gathered} \, .
\fen

\item {\bf (Symmetry nucleation)}  Given a graph, an invertible loop can be nucleated on any face and merged with the bordering edges, where the merging can be regarded as attachments followed by contractions.  For example, on a triangular face,
\[
\begin{gathered}
\begin{tikzpicture}[scale=1]
\draw [line] (-1,0) -- (-1.43,-.25);
\draw [line] (1,0) -- (1.43,-.25);
\draw [line] (0,1.73) -- (0,2.23);
\draw [line,-<-=.53] (-1,0) -- (0,0) node [below] {$\cL_2$} -- (1,0);
\draw [line,->-=.53] (-1,0) -- (-.5,.87) node [above left=-2pt] {$\cL_1$} -- (0,1.73);
\draw [line,-<-=.53] (1,0) -- (.5,.87) node [above right=-2pt] {$\cL_3$} -- (0,1.73);
\draw [line,dashed,-<-=.27] (0,.58) circle (.375);
\node at (0,.58) {$\eta$};
\end{tikzpicture}
\end{gathered}
\quad = \quad
\begin{gathered}
\begin{tikzpicture}[scale=1]
\draw [line] (-1,0) -- (-1.43,-.25);
\draw [line] (1,0) -- (1.43,-.25);
\draw [line] (0,1.73) -- (0,2.23);
\draw [line,-<-=.53] (-1,0) -- (0,0) node [below] {$\eta\cL_2$} -- (1,0);
\draw [line,->-=.53] (-1,0) -- (-.5,.87) node [above left=-2pt] {$\eta\cL_1$} -- (0,1.73);
\draw [line,-<-=.53] (1,0) -- (.5,.87) node [above right=-2pt] {$\eta\cL_3$} -- (0,1.73);
\draw [line,dashed,->-=.6] (.66,.58) -- (.33,0);
\draw [line,dashed,->-=.6] (-.33,0) -- (-.66,.58);
\draw [line,dashed,->-=.6] (-.33,1.07) -- (.33,1.07);
\node at (0,.58) {$\eta$};
\end{tikzpicture}
\end{gathered}
\quad = \quad
\begin{gathered}
\begin{tikzpicture}[scale=1]
\draw [line] (-1,0) -- (-1.43,-.25);
\draw [line] (1,0) -- (1.43,-.25);
\draw [line] (0,1.73) -- (0,2.23);
\draw [line,-<-=.53] (-1,0) -- (0,0) node [below] {$\eta\cL_2$} -- (1,0);
\draw [line,->-=.53] (-1,0) -- (-.5,.87) node [above left=-2pt] {$\eta\cL_1$} -- (0,1.73);
\draw [line,-<-=.53] (1,0) -- (.5,.87) node [above right=-2pt] {$\eta\cL_3$} -- (0,1.73);
\end{tikzpicture}
\end{gathered} \, .
\]
\end{enumerate}

\section{Polynomials with $F$-symbols as roots}
\label{Sec:Poly}

\addtocontents{toc}{\setcounter{tocdepth}{-10}}

\subsection{$G = \bZ_7$}

\ien
P^{\bZ_7}_y(y) =& \, 117649y^{12}-453789y^{11}+1145277y^{10}-1070503y^9+882588y^8-284732y^7 \\ 
& -89977y^6+31488y^5-1828y^4-849y^3+381y^2+45y-1 \, ,
\\
P^{\bZ_7}_z(z) =& \, 343z^6+196z^5-371z^4+27z^3+56z^2-9z-1 \, ,
\\
P^{\bZ_7}_w(w) =& \, 
49w^4-63w^3+15w^2+10w-4 \, .
\fen

\subsection{$G = \bZ_9$}

\ien
P^{\bZ_9}_y(y) =& \ 282429536481y^{24}-2541865828329y^{23}+13891349053584y^{22}-42375665666331y^{21}
\\
&+93048845085738y^{20}-163017616751046y^{19}+191382870385035y^{18}
\\
&-91749046865085y^{17}-71565147070767y^{16}+121393466114850y^{15}
\\
&
-42556511453652y^{14}-23330326470255y^{13}+20787803433577y^{12}
\\
&
-1805958554210y^{11}-2533403044422y^{10}+632950992624y^9
\\
&
+91558817982y^8-30315392921y^7-4655443748y^6+986603649y^5\\
&
+182920180y^4-28268573y^3-1118977y^2-127236y-1801 \, ,
\\
P^{\bZ_9}_z(z) =& \ 
531441z^{12}+885735z^{11}-1535274z^{10}-121014z^9+647352z^8-79407z^7\\
&
-92863z^6+18139z^5+4928z^4-1208z^3-64z+25z-1 \, ,
\\
P^{\bZ_9}_w(w) =& \ 
282429536481 w^{24}-1129718145924 w^{23}+1997927461773 w^{22}-1984755165147 w^{21}
\\
&+1330918519878w^{20}-791614850283 w^{19}+459695402118 w^{18}-222483700269 w^{17}
\\
&+99182263023 w^{16}-47943836820
w^{15}+17026501158 w^{14}-3348784053 w^{13}
\\
&+1374949378 w^{12}-621445880 w^{11}-329500476
w^{10}+412571852 w^9 -148134014 w^8
\\
&+18260969 w^7+2110023 w^6-806198 w^5+47683 w^4+6215 w^3-711
w^2+4 w+1 \, .
\\
P^{\bZ_9}_r(r) =& \ 6561 r^8-8019 r^7+1377 r^6-792 r^5+3349 r^4+4 r^3-662 r^2+52 r+19 \, ,
\\
P^{\bZ_9}_s(s) =& \ 81 s^4+99 s^3+17 s^2-14 s-4 \, .
\fen

\eject

\subsection{$G = \bZ_{11}$}

For the unitary orbit with two solutions
\ien
P^{\bZ_{11}}_{2|y}(y) =& \ 121 y^4-209 y^3+82 y^2+24 y-9 \, ,
\\
P^{\bZ_{11}}_{2|z}(z) =& \ 11 z^2+7 z+1 \, ,
\\
P^{\bZ_{11}}_{2|w}(w) =& \ 121 w^4-88 w^3+38 w^2-13 w+1 \, .
\fen
For the unitary orbit with ten solutions,
\ien
P^{\bZ_{11}}_{10|y}(y) =& \ 25937424601 y^{20}-47158953820 y^{19}+1064291844165 y^{18}+4808654315960 y^{17}
\\
&
+35564388240370 y^{16}+114903432126461 y^{15}+194232171940290 y^{14}
\\
&
+126582540515475 y^{13}-21851286302395 y^{12}-65093840585730y^{11}
\\
&
-20230205549333 y^{10}+6813959963720 y^9+4785911566905 y^8+360322446200 y^7
\\
&
-303249779065 y^6-76228721396 y^5-379548930 y^4+2142467760 y^3
\\
&+324308000 y^2+19299130 y+40207 \, ,
\\
P^{\bZ_{11}}_{10|z}(z) =& \ 
161051 z^{10}+658845 z^9-971630 z^8-542080 z^7+322135 z^6\\
&
+105612 z^5-39815 z^4-6570 z^3+1960 z^2+70 z-19 \, ,
\\
P^{\bZ_{11}}_{10|w}(w) =& \ 25937424601 w^{20}-176846076825 w^{19}+592702305965 w^{18}-1134445659765 w^{17}
\\
&+1534818445765 w^{16}-1765089648718 w^{15}+1769544129045 w^{14}
\\
&-1394768735745 w^{13}+776013578560 w^{12}-263088585485 w^{11}+20179458718
w^{10}
\\
&+32370728245 w^9-20820136235 w^8+6982550700 w^7-1450721110 w^6
\\
&+175316847 w^5-7539540 w^4-877925 w^3+133550 w^2-5960 w+71 \, .
\fen

\subsection{$G = \bZ_{13}$}

\ien
P^{\bZ_{13}}_y(y) =& \ 23298085122481 y^{24}+80647217731665 y^{23}+3069557179509834 y^{22}\\
&
+41919543603471508 y^{21}+536909384312855190 y^{20}+4259352400707950897 y^{19}\\
&
+19179161744641728596 y^{18}+47561155144008593243y^{17}+63626358551986353149 y^{16}
\\
&
+40207662041712799114 y^{15}+1257635216859228766 y^{14}-13522223195096193305 y^{13}\\
&
-6598116247933199625 y^{12}+128413711306511340 y^{11}+938990747292838888y^{10}\\
&
+202797783582401196 y^9-32756778784407789 y^8-16526752437401584 y^7\\
&
-933201395423678 y^6+349378912529867 y^5+53761577382743 y^4+1555890743172 y^3
\\
&
-87453542726 y^2-2773486466 y+28678361 \, ,
\\
P^{\bZ_{13}}_z(z) =& \  
4826809 z^{12}+34901542 z^{11}-124183228 z^{10}-57416398 z^9+51122838 z^8+3476850 z^7
\\
&
-4988283 z^6+418090 z^5+93250 z^4-14139 z^3+205 z^2+38 z-1 \, ,
\\
P^{\bZ_{13}}_w(w) =& \ 
23298085122481 w^{24}-268824059105550 w^{23}+1610738618763716 w^{22}
\\
&
-4730805028787149 w^{21}+8265875258850053 w^{20}-9798763675027379 w^{19}
\\
&
+8948312751528579 w^{18}-6464842564613641 w^{17}+3087209293878385w^{16}
\\
&
-284952516401007 w^{15}-771813881083466 w^{14}+531872957583864 w^{13}
\\
&
-107361616574952 w^{12}-39739582655570 w^{11}+27485052167132 w^{10}
\\
&
-4323332693485 w^9-1159653323459 w^8+583780092624 w^7-51758752951w^6
\\
&
-19939454943 w^5+4746063302 w^4+131285807 w^3-111025779 w^2+2170222 w
\\
&+898159 \, ,
\\
P^{\bZ_{13}}_s(s) =& \ 
28561 s^8-24167 s^7+163930 s^6-225693 s^5+119817 s^4-26999 s^3+1045 s^2+546 s
\\
&-67 \, .
\fen

\bibliography{refs}
\bibliographystyle{JHEP}

\end{document}